\documentclass[a4paper]{amsart}
\usepackage{amsfonts,amsmath,amsthm,amssymb,color}
\usepackage[all]{xy}

\usepackage[ colorlinks=true, linkcolor=blue, citecolor=blue]{hyperref}

\numberwithin{equation}{section}

\begin{document}

\newtheorem{theorem}{Theorem}[section]
\newtheorem{lemma}[theorem]{Lemma}
\newtheorem{proposition}[theorem]{Proposition}
\newtheorem{corollary}[theorem]{Corollary}

\theoremstyle{definition}
\newtheorem{definition}[theorem]{Definition}
\newtheorem{example}[theorem]{Example}

\theoremstyle{remark}
\newtheorem{remark}[theorem]{Remark}
\newtheorem*{ack}{Acknowledgments}

\newenvironment{magarray}[1]
{\renewcommand\arraystretch{#1}}
{\renewcommand\arraystretch{1}}

\newcommand{\mapor}[1]{\smash{\mathop{\longrightarrow}\limits^{#1}}}
\newcommand{\mapin}[1]{\smash{\mathop{\hookrightarrow}\limits^{#1}}}
\newcommand{\mapver}[1]{\Big\downarrow
\rlap{$\vcenter{\hbox{$\scriptstyle#1$}}$}}
\newcommand{\liminv}{\smash{\mathop{\lim}\limits_{\leftarrow}\,}}

\newcommand{\Set}{\mathbf{Set}}
\newcommand{\Art}{\mathbf{Art}}
\newcommand{\solose}{\Rightarrow}

\newcommand{\specif}[2]{\left\{#1\,\left|\, #2\right. \,\right\}}

\renewcommand{\bar}{\overline}
\newcommand{\de}{\partial}
\newcommand{\debar}{{\overline{\partial}}}
\newcommand{\per}{\!\cdot\!}
\newcommand{\Oh}{\mathcal{O}}
\newcommand{\sA}{\mathcal{A}}
\newcommand{\sB}{\mathcal{B}}
\newcommand{\sC}{\mathcal{C}}
\newcommand{\sD}{\mathcal{D}}
\newcommand{\sE}{\mathcal{E}}
\newcommand{\sF}{\mathcal{F}}
\newcommand{\sG}{\mathcal{G}}
\newcommand{\sH}{\mathcal{H}}
\newcommand{\sI}{\mathcal{I}}
\newcommand{\sJ}{\mathcal{J}}
\newcommand{\sL}{\mathcal{L}}
\newcommand{\sM}{\mathcal{M}}
\newcommand{\sP}{\mathcal{P}}
\newcommand{\sU}{\mathcal{U}}
\newcommand{\sV}{\mathcal{V}}
\newcommand{\sX}{\mathcal{X}}
\newcommand{\sY}{\mathcal{Y}}
\newcommand{\sN}{\mathcal{N}}
\newcommand{\sZ}{\mathcal{Z}}

\newcommand{\Aut}{\operatorname{Aut}}
\newcommand{\Mor}{\operatorname{Mor}}
\newcommand{\Def}{\operatorname{Def}}
\newcommand{\Hom}{\operatorname{Hom}}
\newcommand{\Hilb}{\operatorname{Hilb}}
\newcommand{\HOM}{\operatorname{\mathcal H}\!\!om}
\newcommand{\DER}{\operatorname{\mathcal D}\!er}
\newcommand{\Spec}{\operatorname{Spec}}
\newcommand{\Der}{\operatorname{Der}}
\newcommand{\End}{{\operatorname{End}}}
\newcommand{\END}{\operatorname{\mathcal E}\!\!nd}
\newcommand{\Image}{\operatorname{Im}}
\newcommand{\coker}{\operatorname{coker}}
\newcommand{\tot}{\operatorname{tot}}
\newcommand{\Diff}{\operatorname{Diff}}
\newcommand{\ten}{\otimes}
\newcommand{\mA}{\mathfrak{m}_{A}}

\renewcommand{\Hat}[1]{\widehat{#1}}
\newcommand{\dual}{^{\vee}}
\newcommand{\desude}[2]{\dfrac{\de #1}{\de #2}}
\newcommand{\sK}{\mathcal{K}}
\newcommand{\A}{\mathbb{A}}
\newcommand{\N}{\mathbb{N}}
\newcommand{\R}{\mathbb{R}}
\newcommand{\Z}{\mathbb{Z}}
\renewcommand{\H}{\mathbb{H}}
\renewcommand{\L}{\mathbb{L}}
\newcommand{\proj}{\mathbb{P}}
\newcommand{\K}{\mathbb{K}\,}
\newcommand\C{\mathbb{C}}
\newcommand\T{\mathbb{T}}
\newcommand\Del{\operatorname{Del}}
\newcommand\Tot{\operatorname{Tot}}
\newcommand\Grpd{\mbox{\bf Grpd}}
\newcommand\rif{~\ref}

\newcommand\vr{``}
%%%%%%%%%%%%%%%%%%%%%%%%%%%%%%%%%%%%%%%%%%%%%%%%%%%%%%%%%%%%%%%%%%%%5555
\newcommand{\rh}{\rightarrow}
\newcommand{\contr}{{\mspace{1mu}\lrcorner\mspace{1.5mu}}}

\newcommand{\bi}{\boldsymbol{i}}
\newcommand{\bl}{\boldsymbol{l}}

\newcommand{\MC}{\operatorname{MC}}
\newcommand{\Coder}{\operatorname{Coder}}
\newcommand{\CE}{\operatorname{CE}}
\newcommand{\TW}{\operatorname{TW}}
\newcommand{\id}{\operatorname{id}}
\newcommand{\ad}{\operatorname{ad}}

\title{Formality of Kapranov's brackets in K\"ahler geometry via pre-Lie deformation theory}
\author{Ruggero Bandiera}
\address{\newline
Universit\`a degli studi di Roma La Sapienza,\hfill\newline
Dipartimento di Matematica \lq\lq Guido
Castelnuovo\rq\rq,\hfill\newline
P.le Aldo Moro 5,
I-00185 Roma, Italy.}
\email{bandiera@mat.uniroma1.it}

%\date{May 9, 2015.}

\begin{abstract} We recover some recent results by Dotsenko, Shadrin and Vallette on the Deligne groupoid of a pre-Lie algebra, showing that they follow naturally by a pre-Lie variant of the PBW Theorem. As an application, we show that Kapranov's $L_\infty$ algebra structure on the Dolbeault complex of a K\"ahler manifold is homotopy abelian and independent on the choice of K\"ahler metric up to an $L_\infty$ isomorphism, by making the trivializing homotopy and the $L_\infty$ isomorphism explicit.\end{abstract}

\maketitle

\section{Introduction}

Pre-Lie algebras were introduced by Gerstenhaber \cite{Ger} in his study of the deformation theory of associative algebras and  independently by Vinberg \cite{vin}, under the name of left symmetric algebras, in connection with problems in differential geometry:  since then they have been employed in a variety of contexts, ranging from numerical analysis to quantum field theory \cite{chrono,ChLiv,esi,oudom-guin}. As in the recent paper \cite{pLDF}, we are mainly concerned with the role of pre-Lie algebras in deformation theory: we invite the reader to compare the following considerations with the ones in the introduction of loc. cit..

According to a well known principle sponsored by Deligne and others in the 80s \cite{Deligneletter,GoMil1,SchSta}, over a field of characteristic zero every infinitesimal deformation problem is controlled by a dg Lie algebra via the associated deformation functor: this is the functor of Artin rings \cite{Schl} sending a local Artin ring $A$ to the set of Maurer-Cartan elements in $L\ten\mathfrak{m}_A$ (where $\mathfrak{m}_A\subset A$ is the maximal ideal) modulo Gauge equivalence, cf. \cite{GoMil1,iacMan2,KoSo,ManettiRendiconti,BM}.  This principle has been recently made into a rigorous theorem in the papers \cite{LurieDAG,PridUDDT}, cf. also \cite{EDF}: to do so, we have to consider the ordinary deformation functor as the truncation of a derived deformation functor going from the homotopy category of \emph{dg} local Artin rings to the  homotopy category of $\infty$-groupoids, that is, Kan complexes. Explicit models of the derived deformation functor associated to a dg Lie algebra $L$ were introduced by Hinich \cite{hinichdescent,hinichdgC} and Getzler \cite{Getzler04}, following ideas from rational homotopy theory \cite{sullivan}. Hinich's model sends a dg local Artin ring $A$ to the Kan complex of Maurer-Cartan forms on the standard cosimplicial simplex $\Delta_\bullet$ with coefficients in $L\ten\mathfrak{m}_A$. By the proof of the de Rham Theorem given in \cite{dupont}, integration of forms over simplexes is a simplicial quasi-isomorphism from the simplicial dg Lie algebra of forms to the simplicial dg space $C^*(\Delta_\bullet;L\ten\mathfrak{m}_A)$ of non-degenerate \emph{cochains} on $\Delta_\bullet$ with coefficients in $L\ten\mathfrak{m}_A$, hence via homotopy transfer (along Dupont's contraction \cite{dupont,Getzler04}) there is a simplicial $L_\infty$ algebra structure on $C^*(\Delta_\bullet;L\ten\mathfrak{m}_A)$: finally, Getzler's model sends $A$ to the Kan complex $\MC(C^*(\Delta_\bullet;L\ten\mathfrak{m}_A))$ of Maurer-Cartan cochains on $\Delta_\bullet$ with coefficients in $L\ten\mathfrak{m}_A$ (the equivalence of the above definition and the one from \cite{Getzler04} is implicit in loc. cit., an explicit proof can be found in \cite{tesi}). Getzler's model is more closely related to classical deformation theory, for instance, $1$-simplices are in bijective correspondence with arrows in the Deligne groupoid (the action groupoid associated to the Gauge action on Maurer-Cartan elements \cite{GoMil1,hinichdescent,BM}).%, and if $L$ is concentrated in degrees $\geq0$, which is a typical case in deformation theory, it reduces to the nerve of the usual Deligne groupoid . 

We remark here, and hope to elaborate more on this point somewhere else, that the derived deformation theory controlled by a dg associative algebra $B$ should be much simpler to describe than in the general Lie case: like before, as a model of the derived deformation functor we may take $A\to\MC(C^*(\Delta_\bullet;B\ten\mathfrak{m}_A))$, but now we consider $C^*(\Delta_\bullet;B\ten\mathfrak{m}_A)$ as a simplicial dg associative algebra via the cup product. In fact, it should not be hard to show, using the inductive techniques of the recent paper \cite{liemodels}, that this simplicial dg algebra structure is $A_\infty$ isomorphic to the simplicial $A_\infty$ algebra structure induced via homotopy transfer along Dupont's contraction. Using this model, some difficulties in Getzler's theory (for instance, the integration of $\infty$-morphisms or the explicit decription of the $\infty$-groupoid structure in the sense of \cite[Sec. 2]{Getzler04}) may be more easily addressed in the associative case. On the other hand, the class of deformation problems controlled by (the dg Lie algebra associated to) a dg associative algebra is very special, the typical example being the deformations of a complex $(V,d)$, which are controlled by the dg associative algebra $(\End(V),[d,-],\circ)$ (where $\circ$ is the composition product). 

Conversely, for several important deformation problems, such as the deformations of an algebra over an operad \cite{pLDF} or the deformations of the complex structure on a K\"ahler manifold, the Lie bracket of the controlling dg Lie algebra is the commutator of a pre-Lie product. As pre-Lie algebras sit in between associative and Lie algebras, it might be interesting to study the derived deformation theory of pre-Lie algebras along the lines of the associative case sketched above. Of course, the first step in this study, and the only one we will be concerned with in this paper, should be the study of the ordinary deformation theory: this was carried out in \cite[Sec. 4]{pLDF} with different motivations, we refer to loc. cit. for several interesting applications, in particular to homotopy transfer formulas. In Section \ref{sec-prelie} of this work we recover the results from \cite[Sec. 4]{pLDF} using a different method. Whereas in loc. cit. they are proved via a direct computation in a free pre-Lie algebra, where the calculation is done by combinatorics of trees, we shall see that they also follow rather naturally by a pre-Lie variant of the usual Poincar\'e-Birkhoff-Witt Theorem we learned from \cite{esi,oudom-guin}.

More precisely, a pre-Lie algebra structure on $L$ induces a structure of $UL$-bimodule on the symmetric coalgebra $SL$, where we denote by $UL$ the universal enveloping algebra of the associated Lie algebra. The structure of left (resp.: right, according if we are in a left or right pre-Lie algebra) $UL$-module is given by the adjoint of the pre-Lie product, while the structure of right (resp.: left) $UL$-module is given by the associated symmetric braces on $L$ (cf. \cite{oudom-guin}, where there is shown that this induces an isomorphism between the category of pre-Lie algebras and the category of symmetric brace algebras): the usual argument for the PBW Theorem implies that there is an isomorphism of $UL$-bimodules $\eta:UL\to SL$. The morphisms of graded Lie algebras $s,s^\bot:L\to\End(SL)$ associated to the left and right $L$-action respectively factor through the inclusion $\Coder(SL)\to\End(SL)$, which in turn implies that $\eta$ is also an isomorphism of coalgebras. According to the classification theorem from \cite{ChLaz2,Lazarev,methazambon}, the morphism of graded Lie algebras $s-s^\bot:L\times L \to\Coder(SL):(x,y)\to s(x)-s^\bot(y)$ classifies an extension of $L_\infty$ algebras of base $L\times L$ and fiber $L[-1]$. The underlying complex of the total space is naturally isomorphic to the complex $C^\ast(\Delta_1;L)$ of non-degenerate cochains on the $1$-simplex with coefficients in $L$: we denote the $L_\infty$ algebra structure on the total space by $C^\ast(\Delta_1;L)_{p-Lie}$ \footnote{when $L$ is actually an associative algebra, we recover the dg Lie algebra structure on $C^*(\Delta_1;L)_{p-Lie}$ induced by the cup product.} to distinguish it from the previously considered one induced via homotopy transfer along Dupont contraction, which we denote by $C^\ast(\Delta_1;L)_{Lie}$. Working on computations by Fiorenza and Manetti \cite{FMcone}, we show that the latter is classified by an analog morphism of graded Lie algebras $\Phi-\Phi^\bot:L\times L\to\Coder(SL)$, where this time $\Phi,\Phi^\bot:L\to\Coder(SL)$ are associated to the $UL$-bimodule structure on $SL$ induced by the inverse $\operatorname{sym}^{-1}:UL\to SL$ of the usual PBW isomorphism  given by symmetrization. 

We shall denote by $E=\eta\operatorname{sym}$ the composition of the pre-Lie PBW isomorphism $\eta:UL\to SL$ with the usual PBW isomorphism $\operatorname{sym}:SL\to UL$: we call this automorphism of the coalgebra $SL$ the \emph{exponential automorphism}, since it induces the pre-Lie exponential map \cite{chrono,pLDF} on group-like elements. In Theorem \ref{th:deligne} we deduce from the above an explicit isomorphism of $L_\infty$ algebras $C^\ast(\Delta_1;L)_{Lie}\to C^\ast(\Delta_1;L)_{p-Lie}$, closely related to the exponential automorphism $E$. Maurer-Cartan elements in $C^*(\Delta_1;L)_{Lie}$ correspond to arrows in the Deligne groupoid of $L$: in light of the previous considerations, the induced isomorphism $\MC(C^\ast(\Delta_1;L)_{Lie})\to \MC(C^\ast(\Delta_1;L)_{p-Lie})$ of Maurer-Cartan sets recovers the pre-Lie integration of the Deligne groupoid from \cite{pLDF}. In \cite[Sec. 4, Th. 2]{pLDF} the authors give a nice description of the groupoid structure in terms of the symmetric braces on $L$, this can also be recovered in our setting, cf. Theorem \ref{th:formalgrouplaw}: the proof is simplified by the observation (perhaps new) that the morphism of graded Lie algebras $s^\bot:L\to\Coder(SL)$ associated to the symmetric braces is actually a morphism of graded pre-Lie algebras, where $\Coder(SL)$ is a pre-Lie algebra via the Nijenhuis-Richardson product, cf. Proposition \ref{prop:symmbraces}. 

For applications in deformation theory, we are interested in the case when it is given $d$ making $L$ into a dg Lie algebra. We shall denote by $\sK(d)=EdE^{-1}$ the twisting  of the linear coderivation $d:SL\to SL$ by the exponential automorphism $E:SL\to~SL$, and  by $C^\ast(\Delta_1;L,d)_{p-Lie}$ the total space of the $L_\infty$ extension classified by the same morphism $s-s^\bot$ as before, but with fiber the $L_\infty$ algebra structure on $L[-1]$ induced by $\sK(d)$: this fiber identifies with the subalgebra $C^*(\Delta_1,\de\Delta_1; L,d)_{p-Lie}\subset C^*(\Delta_1;L,d)_{p-Lie}$ of cochains relative to the boundary $\de\Delta_1\subset\Delta_1$. Again, we consider the $L_\infty$ algebra $C^\ast(\Delta_1;L,d)_{Lie}$ defined via homotopy transfer along Dupont's contraction: this is the total space of an $L_\infty$ extension classified by $\Phi-\Phi^\bot$, as before, where the fiber is $C^*(\Delta_1,\de\Delta_1;L,d)_{Lie}=(L[-1]-d)$ regarded as an abelian $L_\infty$ algebra (that is, the only non-vanishing bracket is the differential). Finally, the same explicit formulas as in the case $d=0$ define an isomorphsim of $L_\infty$ algebras $C^*(\Delta_1;L,d)_{Lie}\to C^*(\Delta_1;L,d)_{p-Lie}$. In Proposition \ref{prop:main} we show that the higher brackets $\sK(d)_n:L^{\odot n}\to L$, $n\geq1$, may be computed explicitly by a simple recursion \eqref{recursivedefinition}: the latter resembles closely the construction given by Kapranov in the paper \cite{Ka}, and in fact it recovers it as a particular case, cf. the following paragraphs, for this reason we call the brackets $\sK(d)_n$ the Kapranov brackets on $L$ associated to $d$. It might be interesting to point out that, conversely, the brackets associated to $\sK^{-1}(d)=E^{-1}dE$ are a natural generalization to pre-Lie algebras of the classical construction of Koszul brackets on a graded commutative algebra \cite{koszul}, cf. Remark \ref{rem:koszul} and references therein.

In Section \ref{sec-examples} we consider an example from K\"ahler geometry: recall that the deformations of the complex structure on a compact complex manifold $X$ (the prototypical example of a deformation problem)  are controlled by the Kodaira-Spencer dg Lie algebra $(\sA^{0,\ast}(T_X),\debar, [-,-])$ of Dolbeault forms on $X$ with coefficients in the tangent bundle $T_X$ \cite{ManettiRendiconti}. If $X$ is a K\"ahler manifold, the $(1,0)$-component $\nabla$ of the Chern connection (the latter is the only connection compatible with both the metric and the complex structure \cite{Koba}) is flat and torsion-free, in particular, it induces a pre-Lie product on $\sA^{0,\ast}(T_X)$ with associated Lie bracket the usual one. 

In the seminal paper \cite{Ka}, motivated by the study of Rozansky-Witten invariants, Kapranov showed that the Atiyah class makes $T_X[-1]$ into a Lie algebra object in the derived category of bounded below complexes of vector bundles (more in general, coherent sheaves) on $X$: there is a companion theorem for vector bundles, saying that the Atiyah class of a vector bundle $E$ makes $E[-1]$ into a module over the Lie algebra object $T_X[-1]$. He considered two spaces of cochains with coefficients in the (shifted) tangent bundle, and for both choices he showed how the Jacobi identity in the derived category unravels to an actual $L_\infty$ algebra structure on the space of cochains. 

The first construction uses K\"ahler geometry, and the choice of cochains is the usual one of Dolbeault forms: the quadratic bracket is given by the curvature and the higher brackets by its covariant derivatives. We shall see in Theorem \ref{th:main} that the induced $L_\infty$ algebra structure on $\sA^{0,\ast}(T_X)[-1]$ coincides with the one we denoted by $C^*(\Delta_1,\de\Delta_1;\sA^{0,\ast}(T_X),\debar)_{p-Lie}$ in the previous paragraphs. By the results of Section \ref{sec-prelie}, we see in particular that Kapranov's $L_\infty$ algebra is abelian up to homotopy (in fact, it is the loop space of the Kodaira-Spencer algebra in the homotopy category of $L_\infty$ algebras, cf. Remark \ref{rem:loopspace}), and we get moreover an explicit $L_\infty$ isomorphism with the strictly abelian $L_\infty$ algebra $(\sA^{0,\ast}(T_X)[-1],-\debar)=C^*(\Delta_1,\de\Delta_1;\sA^{0,\ast}(T_X),\debar)_{Lie}$ by polarization of the pre-Lie exponential (this may be compared with \cite[Sec. 2.9]{Ka}). As a byproduct of our analysis, we show in Proposition \ref{prop:kapindofmetr} that two different choices of K\"ahler metric induce isomorphic $L_\infty$ algebra structures, by exhibiting an explicit, recursively defined, $L_\infty$ isomorphism.  We should remark that even though Kapranov's brackets on $\sA^{0,\ast}(T_X)[-1]$ are linear over the Dolbeault algebra $\sA_X^{0,\ast}$, and they are independent on the metric up to an $\sA_X^{0,\ast}$-linear $L_\infty$ isomorphism, homotopy abelianity only holds, in general, over the field of complex numbers. 

As pointed out by the referee, both results are expected, and in a certain measure already implicit in Kapranov's second construction, which uses formal geometry  \cite{bott,GF} and as cochains the relative forms with cofficients in the space of formal exponential maps, cf. \cite[Sec. 4]{Ka} (thus, the only actual novelties are the explicit isomorphisms and the method of proof, using pre-Lie algebras) . In particular, the second construction uses only the complex structure, which should imply Proposition \ref{prop:kapindofmetr}. Following a point of view further developed in the papers \cite{Caldararu1,Caldararu2}, Kapranov's construction may be considered as a first step in establishing a dictionary between Lie theory and complex (or algebraic) geometry. In this dictionary, the manifold $X$ corresponds to the Lie object $T_X[-1]$ and the derived category of bounded below complexes of coherent sheaves on $X$ to the category of representations of $T_X[-1]$. Moreover, the structure sheaf $\Oh_X$ corresponds to the trivial representation and the (shifted) tangent sheaf to the adjoint representation: finally, the complex $\sA^{0,\ast}(T_X)[-1]$, regarded as $\mathbf{R}\operatorname{Hom}(\Oh_X,T_X[-1])$ (where we compute $\mathbf{R}\operatorname{Hom}(-,-)$ in the \emph{dg} category of coherent sheaves), corresponds to the Chevalley-Eilenberg complex of $T_X[-1]$ with coefficients in the adjoint representation (cf. \cite{cartan,weibel}), or in other words, to the derived center of the Lie algebra object $T_X[-1]$, and in particular its $L_\infty$ algebra structure should be abelian up to homotopy (over the base field, and not in general over the Chevalley-Eilnberg algebra with coefficients in the trivial representation). Kapranov's constructions have been recently the object of extensive study, for further readings we refer to \cite{Caldararu1,Caldararu2,Xuatiyah,Costello,Hennion,Xuexponential,Xukapdg,Yu}

\begin{ack} It is a pleasure to thank my Ph.D. advisor Marco Manetti for supporting me throughout the writing of this paper, as well as Domenico Fiorenza, Jim Stasheff and the referee for numerous remarks, corrections and suggestions. \end{ack}

\section{Preliminaries on $L_\infty[1]$ algebras}\label{Sec-L_inftyalg}

In the first part of this section, mostly with the aim to fix notations, we recall some well known fact on $L_\infty$ algebras: at the end we recall the classification of $L_\infty$ algebra extensions from \cite{ChLaz2,Lazarev,methazambon}. We work over a field $\K$ of characteristic zero, graded means $\Z$-graded. For a graded space $V=\oplus_{i\in\Z}V^i$ we denote by $V^{\otimes n}$ the $n$-th tensor power of $V$, i.e., the tensor product of $n$ copies of $V$, and by $V^{\odot n}$, resp. $V^{\wedge n}$, the $n$-th symmetric power, resp. exterior power, of $V$, that is, the space of coinvariants of $V^{\otimes n}$ under the natural, resp. alternate, action of the symmetric group $S_n$ (with the usual Koszul rule for twisting signs); $V^{\ten0}=V^{\odot0}=V^{\wedge0}:=\K$. Given an integer $k\in\Z$ we denote by $V[k]$ the shifted space $V[k]^i=V^{k+i}$. For graded spaces $V$ and $W$, we denote by $\Hom(V,W)=\oplus_{i\in\Z}\Hom^i(V,W)$ the internal mapping space in the category of graded spaces. Given integer $i_1+\cdots+i_k=n$ we denote by $S(i_1,\ldots,i_k)\subset S_n$ the set of $(i_1,\ldots,i_k)$-unshuffles, that is, permutations $\sigma\in S_n$ such that $\sigma(1)<\cdots<\sigma(i_1)$, $\sigma(i_1+1)<\cdots<\sigma(i_1+i_2)$, ..., $\sigma(i_1+\cdots+i_{k-1}+1)<\cdots<\sigma(n)$.

Let $V$ be a graded space, we denote by $SV=\oplus_{n\ge0}V^{\odot n}$ (or sometimes $S(V)$) the symmetric coalgebra over~$V$: the graded cocommutative coalgebra structure is given by the unshuffle coproduct $\Delta:SV\rh SV\otimes SV$
\[ \Delta(v_1\odot\cdots\odot v_n)=\sum_{i=0}^{n}\sum_{\sigma\in S(i,n-i)}\varepsilon(\sigma)(v_{\sigma(1)}\odot\cdots\odot v_{\sigma(i)})\otimes(v_{\sigma(i+1)}\odot\cdots\odot v_{\sigma(n)}), \]
where $\varepsilon(\sigma)=\varepsilon(\sigma;v_1,\ldots,v_n)$ is the Koszul sign, and with the understanding that for $i=0$ or $i=n$ we put $1\in\K=V^{\odot0}\subset SV$ in place of the empty string (in particular $\Delta(1)=1\otimes1$). The natural projection $SV\rh V^{\odot0}=\K$ and the natural inclusion $\K=V^{\odot0}\rh SV$ give respectively a counit and a coaugmentation for the coalgebra structure. The reduced symmetric coalgebra $(\overline{SV},\overline{\Delta})$ over $V$ (sometimes denoted by $\overline{S}(V)$) is the space $\oplus_{n\geq1}V^{\odot n}\,=:\,\overline{SV}\subset SV$, with the reduced coproduct $\overline{\Delta}$ defined by the same formula as before but taking the sum from $i=1$ to $n-1$.

Given graded spaces $V,W$ and a morphism of graded coaugmented coalgebras $F:SV\rh SW$, the Taylor coefficients of $F$ are the morphism of graded spaces $f_n:V^{\odot n}\hookrightarrow SV\xrightarrow{F}SW\xrightarrow{p}W$, $n\geq1$, where we denote by $p:SW\rh W$ the natural projection. We can reconstruct $F$ knowing its Taylor coefficients via $F(1)=1$ and for $n\geq1$
\begin{equation}\label{invtocorestriction} F(v_1\odot\cdots\odot v_n)=\sum_{k=1}^{n}\frac{1}{k!}\sum_{i_1+\cdots +i_k=n}\sum_{\sigma\in S(i_1,\ldots,i_k)}\varepsilon(\sigma)f_{i_1}(v_{\sigma(1)}\odot\cdots)\odot \cdots\odot f_{i_k}(\cdots\odot v_{\sigma(n)}),  \end{equation}
conversely, any family of morphisms $f_n:V^{\odot n}\rh W$, $n\ge1$, determines a morphism of graded coaugmented coalgebras $F$ in this way. Notice that $F(\overline{SV})\subset\overline{SW}$.

 Similarly, every coderivation $Q:SV\rh SV$ is determined by its corestriction $SV\xrightarrow{Q}SV\xrightarrow{p}V$, giving an isomorphism of graded spaces
\[ \Coder(SV)\xrightarrow{\cong} \Hom(SV,V)=\prod_{n\geq0}\Hom(V^{\odot n},V):Q\xrightarrow{\qquad}pQ=(q_0,q_1,\ldots,q_n,\ldots).  \]
The linear maps $q_n:V^{\odot n}\rh V$, $n\geq0$, are called the Taylor coefficients of $Q$, we reconstruct $Q$ from its Taylor coefficients via
\begin{equation}\label{inversetocorestriction}Q(v_1\odot\cdots\odot v_n)=\sum_{i=0}^n\sum_{\sigma\in S(i,n-i)}\varepsilon(\sigma)q_{i}(v_{\sigma(1)}\odot\cdots\odot v_{\sigma(i)})\odot v_{\sigma(i+1)}\odot\cdots\odot v_{\sigma(n)},\end{equation}
always with the understanding that $q_0(\emptyset):=q_0(1)$ (in particular $Q(1)=q_0(1)\in V\subset SV$).
\begin{remark}\label{costcoder}We call a coderivation $Q\in\Coder(SV)$ \emph{linear} (resp.: \emph{constant}) if $q_n=0$ for $n\neq1$ (resp.: $n\neq0$). Given $v\in V$, we denote by $\Coder(SV)\ni \sigma_v=v\odot-:SV\to SV$ the constant coderivation with $\sigma_v(1)=v$.
\end{remark}
Similarly, a coderivation $Q:\overline{SV}\to\overline{SV}$ is determined by its corestriction $pQ:\overline{SV}\to V$, and there is  an embedding $\Coder(\overline{SV})\to\Coder(SV)$ with image the graded Lie subalgebra of coderivations with vanishing constant term.

The graded Lie algebra structure on $\Coder(SV)$ is induced by a right pre-Lie product (cf. Definition \ref{def-prelie}), called the \emph{Nijenhuis-Richardson product}, which we denote by $\circ$: it sends $Q,R\in\Coder(SV)$ to the only coderivation $Q\circ R$ which corestricts to $pQR:SV\rh V$. In Taylor coefficients, if $pQ=(q_0,\ldots,q_n,\ldots)$ and $pR=(r_0,\ldots,r_n,\ldots)$, then \begin{equation}\label{Nij.-Ric.product} (Q\circ R)_n(v_1\odot\cdots\odot v_n) = \sum_{i=0}^{n}\sum_{\sigma\in S(i,n-i)}\varepsilon(\sigma) q_{n-i+1}(r_i(v_{\sigma(1)}\odot\cdots\odot v_{\sigma(i)})\odot v_{\sigma(i+1)}\odot\cdots\odot v_{\sigma(n)}). \end{equation}
As already remarked, the associated Lie bracket, called the \emph{Nijenhuis-Richardson bracket}, is the usual commutator of coderivations.
\begin{remark}\label{brwithcostcoder} Given $Q\in\Coder(SV)$ and a constant coderivation $\sigma_v$, the Nijenhuis-Richardson bracket $[Q,\sigma_v]=Q\circ\sigma_v-(-1)^{|Q||v|}\sigma_v\circ Q=Q\circ\sigma_v$ is given in Taylor coefficients by $[Q,\sigma_v]_n=q_{n+1}(v\odot-)$: more precisely,
\[ [Q,\sigma_v]_n(v_1\odot\cdots\odot v_n)=q_{n+1}(v\odot v_1\odot\cdots\odot v_n)\quad\mbox{for all $n\geq0$}.\] 
In particular, the constant coderivations span an abelian Lie subalgebra of $\Coder(SL)$. Given $f:V\rh V$, regarded as a linear coderivation on $SV$, we find $[f,\sigma_v]=\sigma_{f(v)}$. 
\end{remark}

\begin{definition} An $L_\infty[1]$ algebra structure on a graded space $V$ is the datum of a dg coalgebra structure on the reduced symmetric coalgebra $\overline{SV}$: in other words, is the datum of a coderivation $Q=(0,q_1,\ldots,q_n,\ldots)$, $|Q|=1$, such that $Q\circ Q=0$. A morphism $F:(V,Q)\rh(W,R)$ of $L_\infty[1]$ algebras is a morphism $F=(f_1,\ldots,f_n,\ldots):(SV,Q)\rh(SW,R)$ of dg coaugmented coalgebras: it is \emph{strict} if $f_n=0$ for $n\geq2$.
\end{definition}

\begin{remark}\label{decalalge} We recall the link with $L_\infty$ algebras, as defined for instance in \cite{LaSt}. We denote by $s^{-1}:V\rh V[1]$ and $s:V[1]\rh V$ the shifts, then d\'ecalage
\[\mbox{d\'ec}:\Hom(V^{\ten n},W)\rh\Hom(V[1]^{\ten n},W[1]):f\rh s^{-1}\, f\, s^{\ten n}\] is an isomorphism of vector spaces which shifts the degrees by $n-1$. Taking into account the signs coming from the Koszul rule, it can be checked that it restricts to a degree $n-1$ isomorphism $\mbox{d\'ec}:\Hom(V^{\wedge n},W)\rh\Hom(V[1]^{\odot n},W[1])$. An $L_\infty$ algebra structure on $V$ is the datum of a family of degree $2-n$ graded antisymmetric bracket $l_n:V^{\wedge n}\rh V$, $n\geq1$, satisfying some relations: these translate into the requirement that the $q_n:=\mbox{d\'ec}(l_n):V[1]^{\odot n}\rh V[1]$ are the Taylor coefficients of an $L_\infty[1]$ algebra structure on $V[1]$. Similarly, an $L_\infty$ morphism between $L_\infty$ algebras $V$ and $W$ is a family of degree $1-n$ maps $f_n:V^{\wedge n}\rh W$, $n\geq1$, such that $\mbox{d\'ec}(f_n):V[1]^{\odot n}\rh W[1]$ are the Taylor coefficients of an $L_\infty[1]$ morphism between the corresponding $L_\infty[1]$ algebras. In other words, d\'ecalage is an isomorphism between the categories of $L_\infty$ and $L_\infty[1]$ algebras, cf. \cite{LaSt} for more details.\end{remark}

\begin{example}\label{eq:quilconstr}
In particular, a dg Lie algebra structure $(L,d,[\cdot,\cdot])$ on $L$ induces an $L_\infty[1]$ algebra structure $Q$ on $L[1]$ with Taylor coefficients $q_1(s^{-1}x)=-s^{-1}dx$, $q_2(s^{-1}x\odot s^{-1}y)=(-1)^{|x|}s^{-1}[x,y]$ and $q_n=0$ for $n\geq 3$.
\end{example}

\begin{definition} Given an $L_\infty[1]$ algebra structure $Q$ on $V$, the linear Taylor coefficient $q_1$ satisfies $q_1^2=0$: we call the dg space $(V,q_1)$ the \emph{tangent complex} of the $L_\infty[1]$ algebra $(V,Q)$, and its cohomology $H(V,q_1)$ the \emph{tangent cohomology} of $(V,Q)$.  Given $F:(V,Q)\rh(W,R)$ an $L_\infty[1]$ morphism, its linear part is a dg morphism between the tangent complexes $f_1:(V,q_1)\rh(W,r_1)$: if $f_1$ is a quasi-isomorphism then $F$ is called a \emph{weak equivalence}. An $L_\infty[1]$ algebra $(V,Q)$ is called \emph{abelian} if $Q$ is a linear coderivation, and is called \emph{homotopy abelian} if it is weakly equivalent to an abelian $L_\infty[1]$ algebra.
\end{definition}

\begin{remark} A minimal model of the $L_\infty[1]$ algebra $(V,Q)$ is the datum of an $L_\infty[1]$ weak equivalence $F:(W,R)\rh(V,Q)$ with $(W,R)$ minimal: recall that this means $r_1=0$. Structure theory of $L_\infty[1]$ algebras (cf. \cite{KoSo}) says that a minimal model of $(V,Q)$ always exists and it is well defined up to a non canonical $L_\infty[1]$ isomorphism over $(V,Q)$, moreover, $(V,Q)$ is isomorphic, as an $L_\infty[1]$ algebra, to the direct product of a minimal model and an acyclic complex, the latter regarded as an abelian $L_\infty[1]$ algebra. It is easy to show that $(V,Q)$ is homotopy abelian if and only if the $L_\infty[1]$ algebra structure on a minimal model is trivial if and only if the following seemingly stronger condition holds: there is an $L_\infty[1]$ isomorphism $F:(V,q_1)\rh(V,Q)$ with $f_1=\id_V$, where in the left hand side we regard $(V,q_1)$ as an abelian $L_\infty[1]$ algebra. We refer the reader to \cite{Manettiformality} for a more exhaustive discussion on homotopy abelian $L_\infty[1]$ algebras and their role in deformation theory.
\end{remark}

Given an $L_\infty[1]$ algebra $(V,Q)$, under suitable conditions ensuring convergence, for instance, if the Taylor coefficients of $Q$ are continuous with respect to a complete descending filtration $V=F^1V\supset\cdots\supset F^pV\supset\cdots$ on $V=\,\underleftarrow{\operatorname{lim}}\,V/F^pV$, it makes sense to consider the set $\MC(V)$ of solutions of the \emph{Maurer-Cartan equation}
\[\MC(V)\,:=\,\{ x\in V^0 \,\,\mbox{s.t.}\,\,\sum_{n\geq1}\frac{1}{n!}q_n(x\odot \cdots\odot x) = 0 \}\]
\begin{remark}\label{rem:Gauge} Given a dg Lie algebra structure $(L,d,[-,-])$ on $L$, seen as an $L_\infty[1]$ algebra structure on $L[1]$, the Maurer-Cartan equation is the usual one $\MC(L)=\{ x\in L^1\,\,\mbox{s.t.}\,\,dx+\frac{1}{2}[x,x]=0 \}$. As well known, the exponential group of $L$, namely, $L^0$ equipped with the group structure given by the Baker-Campbell-Hausdorff product $\bullet$, acts on the set $\MC(L)$ via the \emph{Gauge action} \cite{GoMil1}: we shall denote the latter by $e^-\ast_G -:L^0\times\MC(L)\to \MC(L):(a,x)\to e^a\ast_G x$, explicitly, where $\ad_a=[a,-]:L\to L$ is the adjoint,
\[e^a\ast_G x = x+\sum_{n\geq0}\frac{(\ad_a)^n}{(n+1)!}([a,x]-da). \] 
The \emph{Deligne groupoid} of $L$, of fundamental importance in the study of deformation theory via dg Lie algebras \cite{GoMil1,hinichdescent,iacMan2,BM}, is the action groupoid associated to the Gauge action: the objects are the Maurer-Cartan elements $x\in\MC(L)$ and the arrows the Gauge equivalences between them, the composition is given by the Baker-Campbell-Hausdorff product.
\end{remark}

We close this section by reviewing the classification of $L_\infty$ extensions from \cite{ChLaz2,Lazarev,methazambon}.
\begin{definition}
A graded subspace $I\subset L$ of the $L_\infty[1]$ algebra $(V,q_1,\ldots,q_k,\ldots)$ is an $L_\infty[1]$ ideal if $q_k(I\ten V^{\odot k-1})\subset I$ for all $k\geq1$: then there is an induced $L_\infty[1]$ algebra structure on $V/I$ such that the projection $V\rh V/I$ is a strict morphism of $L_\infty[1]$ algebras. The sequence of $L_\infty[1]$ algebras and strict morphisms $I\rh V\rh V/I$ is called an \emph{$L_\infty[1]$ extension} of fiber $I$ and base $V/I$.
\end{definition}

\begin{example}\label{ex:extensions}
Let $(W,r_1,\ldots,r_n,\ldots)$, $(I,q_1,\ldots,q_n,\ldots)$ be $L_\infty[1]$ algebras. We consider the dg Lie algebra structure on $\Coder(SI)$ given by the Nijenhuis-Richardson bracket and the differential $[Q,-]$, together with the associated $L_\infty[1]$ algebra structure on $\Coder(SI)[1]$. Given a morphism $F=(f_1,\ldots,f_k,\ldots):W\rh\Coder(SI)[1]$ of $L_\infty[1]$ algebras, this defines an $L_\infty[1]$ extension $I\rh W\times_F I\rh W$ of fiber $I$ and base $W$ as follows: the underlying space is $W\times_F I=W\times I$, the $L_\infty[1]$ structure is given in Taylor coefficients $\widetilde{q}_n:(W\times I)^{\odot n}\rh W\times I$ by (according to the decomposition $\Hom((W\times I)^{\odot n},W\times I)=\prod_{j+k=n}\Hom(W^{\odot j}\ten I^{\odot k},W\times I)$)
\[\widetilde{q}_k(i_1\odot\cdots\odot i_k)=(0,q_k(i_1\odot\cdots\odot i_k)),\]
\[ \widetilde{q}_j(w_1\odot\cdots\odot w_j)=(r_j(w_1\odot\cdots\odot w_j), sf_j(w_1\odot\cdots\odot w_j)_0(1)),\]
\[ \widetilde{q}_{j+k}(w_1\odot\cdots\odot w_j\ten i_1\odot\cdots\odot i_k )=(0,sf_j(w_1\odot\cdots\odot w_j)_k(i_1\odot\cdots\odot i_k)),\]
where $sf_j$ is the composition $W^{\odot j}\xrightarrow{f_j}\Coder(SI)[1]\xrightarrow{s}\Coder(SI)$. We refer to \cite{ChLaz2,Lazarev,methazambon} for a proof that $\widetilde{Q}=(\widetilde{q}_1,\ldots,\widetilde{q}_n,\ldots)$ is an $L_\infty[1]$ algebra structure, then it's clear that  $I\rh W\times_F I\rh W$ is an $L_\infty[1]$ extension of fiber $I$ and base $W$, which we call the $L_\infty[1]$ extension classified by $F$.
\end{example}

Conversely, all $L_\infty[1]$ extensions can be constructed as in the previous example.

\begin{theorem}\label{th:extensions} Given an $L_\infty[1]$ extensions $I\rh L\rh L/I$ and an isomorphism of graded spaces $\varphi:L\to L/I\times I$, there is an $L_\infty[1]$ morphism $F:L/I\rh\CE(I)[1]$ such that $\varphi$ is a strict isomorphism of $L_\infty[1]$ algebras $\varphi:L\to L/I\times_F I$.
\end{theorem}
\begin{proof} Cf. \cite{ChLaz2,Lazarev,methazambon}.\end{proof}

We shall need the following lemma, whose proof is a tedious but direct verification. Given an $L_\infty[1]$ extension $I\to W\times _F I\to W$, classified by $F:W\to\Coder(SI)[1]$ as in the previous example, together with an $L_\infty[1]$ \emph{iso}morphism $G=(g_1,\ldots,g_n,\ldots):(I,Q)\to(I', Q')$, we shall denote by $G_\ast F:W\to\Coder(SI')[1]$ the composition of $F$ and the isomorphism of dg Lie algebras $G-G^{-1}:\Coder(SI)\to\Coder(SI')$, and by $\widetilde{G}:W\times_F I\to W\times_{G_\ast F} I'$ given in Taylor coefficients $\widetilde{g}_n:(W\times I)^{\odot n}\to W\times I'$ by 
\begin{equation}\label{eq:Etilde} \widetilde{g}_1(w,i)=(w,g_1(i)),\qquad\widetilde{g}_n((w_1,i_1)\odot\cdots\odot(w_n,i_n))=(0,g_n(i_1\odot\cdots\odot i_n))\,\,\,\mbox{for all $n\geq2$}.\end{equation}
\begin{lemma}\label{lem:extensions} The diagram
\[ \xymatrix{ I\ar[r]\ar[d]_G & W\times_F I\ar[r]\ar[d]_{\widetilde{G}} & W\ar@{=}[d] \\ I'\ar[r] & W\times_{G_* F} I'\ar[r] & W } \]
is an isomorphism of $L_\infty$ extensions.
\end{lemma}

\section{Pre-Lie deformation theory and PBW Theorem}\label{sec-prelie}

%Pre-Lie algebras were introduced by Gerstenhaber \cite{Ger} in his study of the deformation theory of associative algebras and  independently by Vinberg \cite{vin}, under the name of left symmetric algebras, in connection with problems in differential geometry.

\begin{definition}\label{def-prelie} A \emph{graded left pre-Lie algebra} $(L,\triangleright)$ is a graded space $L$ together with a bilinear product $\triangleright:L^{\ten 2}\rh L$ such that the associator, defined by
\[A:L^{\ten3}\rh L:x\ten y\ten z\rh A(x,y,z) = (x\triangleright y)\triangleright z - x\triangleright(y\triangleright z),\]
is graded symmetric in the first two arguments, that is,
\begin{equation}\label{preLiecondition1}A(x,y,z)=(-1)^{|x||y|}A(y,x,z),\qquad\forall x,y,z\in L.\end{equation}
As well known, this implies that the commutator
\[ [\cdot,\cdot]:L^{\wedge 2}\rh L,\qquad[x,y]:=x\triangleright y- (-1)^{|x||y|} y\triangleright x,\]
satisfies the graded Jacobi identity, hence defines a graded Lie algebra structure on $L$. Motivated by geometric examples, we denote by $\nabla:L\rh\End(L):x\rh\{\nabla_x:y\rh x\triangleright y\}$ the left adjoint morphism; then the left pre-Lie identity \eqref{preLiecondition1} is equivalent to 
\begin{equation}\label{preLiecondition2}  [\nabla_x,\nabla_y]=\nabla_{[x,y]}\qquad\mbox{$\forall x,y\in L$}.\end{equation}
A graded space $L$ equipped with a bilinear product $\triangleleft$ such that the associator is graded symmetric in the \emph{last} two variables is called a \emph{graded right pre-Lie algebra}: again, the associated commutator is a Lie bracket. %An important example of right pre-Lie algebra is $\Coder(SL)$ with the Nijenhuis-Richardson product. 
It is straightforwad that $\triangleright$ is a left pre-Lie product on $L$ if and only if the opposite product  $x\triangleleft y:= (-1)^{|x||y|}y\triangleright x$ is right pre-Lie: in the following we shall restrict our attention to left pre-Lie algebras, keeping in mind that all the results can be translated in the right pre-Lie case via the previous observation.
\end{definition}

We denote by $(UL,\ast,\Delta_{UL})$ the universal enveloping algebra of the graded Lie algebra $(L,[\cdot,\cdot])$, with its structure of biaugmented graded cocommutative bialgebra, and as usual by $(SL,\Delta)$ the symmetric coalgebra over $L$. 
\begin{theorem}[Poincar\'e-Birkhoff-Witt]\label{PBW} The symmetrization map $\operatorname{sym}:SL\rh UL$
\begin{equation}\label{sym} \operatorname{sym}(1)=1,\qquad\operatorname{sym}(x_1\odot\cdots\odot x_n)=\frac{1}{n!}\sum_{\sigma\in S_n}\varepsilon(\sigma)x_{\sigma(1)}\ast\cdots\ast x_{\sigma(n)},\end{equation}
is an isomorphism of coaugmented coalgebras.
\end{theorem}
Following the papers \cite{esi,oudom-guin}, we shall recall in the following Theorem \ref{preLiePBW} a pre-Lie variant of the above well known theorem, which will be the main tool behind the computations of this section.

Given a graded left pre-Lie algebra $(L,\triangleright)$, there is an induced Lie action of $L$ on $SL$ by coderivations
\begin{equation}\label{s} s:L\rh\Coder(SL):x\to s(x):=\sigma_x+\nabla_x,\end{equation}
where $\sigma_x=x\odot-:SL\to SL$ is the constant coderivation as in Remark \ref{costcoder}, and we regard $\nabla_x$ as a linear coderivation. This is in fact a morphism of graded Lie algebras: by \eqref{preLiecondition2} and Remark \ref{brwithcostcoder}
\[[s(x),s(y)]=[\sigma_x+\nabla_x,\sigma_y+\nabla_y]=\sigma_{\nabla_x(y)-(-1)^{|x||y|}\nabla_y(x)}+[\nabla_x,\nabla_y]=\sigma_{[x,y]}+\nabla_{[x,y]}=s([x,y]). \]
We can regard $SL$ as a left $UL$-module via the unique extension of \eqref{s} to a morphism of graded algebras $s:UL\rh\End(SL)$. 
\begin{theorem}\label{preLiePBW} The linear map $\eta:UL\rh SL$ 
\begin{equation}\label{eta} \eta(1)=1,\qquad\eta(x_1\ast\cdots\ast x_n)=s(x_1\ast\cdots\ast x_n)(1)=s(x_1)\cdots s(x_n)(1),\end{equation}
is both an isomorphism of coaugmented coalgebras and of left $UL$-modules.
\end{theorem}
\begin{proof} $\eta$ is by construction a morphism of left $UL$-modules: in fact, it is the only morphism of left $UL$-modules such that $\eta(1)=1$. Denoting by $UL_{\leq n}$ and $SL_{\leq n}$ the subspaces spanned by words of lenght $\leq n$, we see inductively that $\eta$ restricts to an isomorphism $UL_{\leq n}\xrightarrow{\cong}SL_{\leq n}$ for all $n\geq0$: in fact, $\eta$ restricts to the identity on $UL_{\leq 1}=\K1\oplus L=SL_{\leq 1}$, while the inductive step follows since $\eta(x_1\ast\cdots\ast x_n)=x_1\odot\cdots\odot x_n+\left\{\mbox{terms in $SL_{\leq n-1}$}\right\}$. $\eta$ is clearly compatible with the counits and the coaugmentations, hence it only remains to show that it is compatible with the coproducts: this follows since $s(x)$ is a coderivation for all $x\in L$, thus
\begin{multline*}
\Delta\eta(x_1\ast\cdots\ast x_n)=\Delta s(x_1)\cdots s(x_n)(1)=\\=\left(s(x_1)\otimes\operatorname{id}+\operatorname{id}\otimes s(x_1)\right)\cdots\left(s(x_n)\otimes\operatorname{id}+\operatorname{id}\otimes s(x_n)\right)(1\otimes 1)=\\ =\sum_{i=0}^n\sum_{\sigma\in S(i,n-i)}\varepsilon(\sigma) s(x_{\sigma(1)}) \cdots s(x_{\sigma(i)})(1) \otimes s(x_{\sigma(i+1)}) \cdots s(x_{\sigma(n)})(1)
=\\=\sum_{i=0}^n\sum_{\sigma\in S(i,n-i)}\varepsilon(\sigma)\eta(x_{\sigma(1)}\ast\cdots\ast x_{\sigma(i)})\ten\eta(x_{\sigma(i+1)}\ast\cdots\ast x_{\sigma(n)})= \\=(\eta\otimes\eta)\left((x_1\otimes 1+1\otimes x_1)\ast\cdots\ast(x_n\otimes 1+1\otimes x_n)\right)=\\=(\eta\otimes\eta)\left(\Delta_{UL}(x_1)\ast\cdots\ast\Delta_{UL}(x_n)\right)=(\eta\otimes\eta)\Delta_{UL}(x_1\ast\cdots\ast x_n).
\end{multline*}
\end{proof}
\begin{remark}\label{rem:PBWdata} The same proof shows the following more general fact: given a graded space $L$ and a linear embedding $s:L\to\Coder(SL)$ such that the image is closed under the Nijenhuis-Richardson bracket and $s(x)_0(1)=x$ for all $x\in L$, then $L$ carries together with the induced graded Lie algebra structure an \vr exotic'' PBW isomorphism $\eta:UL\to SL$, defined as in the claim of the previous proposition. For instance, given a graded post-Lie algebra $(L,[-,-],\triangleright)$, as defined for instance in \cite{postlie}, it is not hard to show (using \cite[Th. 1.2]{HDB}) that $s:x\to\Phi(x)+\nabla_x$ satisfies the above hypotheses, where $\nabla_x=x\triangleright-:L\to L$ and $\Phi$ is induced from the bracket $[-,-]$ as in \eqref{Phi}. This recovers the \vr exotic'' post-Lie PBW isomorphism from \cite{postlie}.
\end{remark}

Given $x\in L$, the operators $x\ast-,-\ast x:UL\to UL$ of left and right multiplication with $x$ are coderivations of the bialgebra $UL$. By construction, we have $s(x)=\eta(x\ast-)\eta^{-1}$ for all $x\in L$, where $\eta$ is the isomorphism from the previous theorem. We shall denote by $s^\bot$ the correspondence $s^\bot:L\to\Coder(SL):x\to s^\bot(x):=\eta(-\ast x)\eta^{-1}$. Since $(UL,\ast)$ is an associative algebra $[-\ast x,y\ast -]=0$ for all $x,y\in L$, which implies that the coderivation $s^\bot(x)$ satisfies $[s^\bot(x),s(y)]=0$ for all $y\in L$, moreover, $s^\bot(x)(1)=\eta(-\ast x)\eta^{-1}(1)=\eta(x)=x$. Conversely, these two properties characterize $s^\bot(x)$ completely, as they determine the Taylor coefficients $s^\bot(x)_n$ recursively (cf. Remark \ref{brwithcostcoder}): for all $n\geq1$ and $x,y\in L$
\begin{multline*}s^\bot(x)_{n+1}(y\odot-)=[s^\bot(x)_{n+1},\sigma_y]=[s^\bot(x)_{n+1},s(y)_0]=\\=[s^\bot(x),s(y)]_n-[s^\bot(x)_n,s(y)_1]=-[s^\bot(x)_n,\nabla_y].\end{multline*}
We put $\{y_1,\ldots, y_n;x\}\,:=\, (-1)^{|x|(|y_1|+\cdots+|y_n|)}s^\bot(x)_n(y_1\odot\cdots\odot y_n)$. When we make the previous recursion explicit, we find
\[\{1;x\}=x,\quad \{y;x\}=y\triangleright x,\quad \{y, z;x\}= y\triangleright (z\triangleright x)-(y\triangleright z)\triangleright x,\]
\[\{y, y_1,\ldots, y_n;x\}=\{y;\{y_1,\ldots, y_n;x\}\} - \sum_{i=1}^n(-1)^{|y|(|y_1|+\cdots+|y_{i-1}|)}\{y_1,\ldots,\{y;y_i \},\ldots,y_n;x\}.\]
In other words, the higher braces $\{-;x\}:SL\to L$ are the ones defining the associated structure of (right) symmetric brace algebra on the (left) pre-lie algebra $L$, cf. \cite{oudom-guin,pLDF}: this follows by a direct comparison with the definitions in loc. cit.. Recall that $\Coder(SL)$ is a right pre-Lie algebra via the Nijenhuis-Richardson product $\circ$, moreover, we may regard the left pre-Lie algebra $(L,\triangleright)$ as a right pre-Lie algebra $(L,\triangleleft)$ via the opposite product $x\triangleleft y\,:=(-1)^{|x||y|}y\triangleright x$.
\begin{proposition}\label{prop:symmbraces} The correspondence $s^\bot:(L,\triangleleft)\to(\Coder(SL),\circ)$ is a morphsim of graded right pre-Lie algebras. For a fixed $x\in L$, the coderivation $s^\bot(x)$ is uniquely determined by the properties
\[ s^\bot(x)_0(1)=x,\qquad[s^\bot(x),s(y)]=0\,\,\mbox{for all $y\in L$.}\]
\end{proposition}
\begin{proof} We have already proved the second statement, we will show that this also implies the first one. Given a coderivation $Q\in\Coder(SL)$, we see by a straightforward computation that the failure of $[Q,-]$ to be a derivation with respect to $\circ$ is measured by the associator $A(Q,-,-)$: more precisely,
\[ A(Q,R,S)=(Q\circ R)\circ S-Q\circ(R\circ S)=[Q,R]\circ S+(-1)^{|Q||R|}R\circ[Q,S]-[Q,R\circ S].\]
Since $A(Q,-,-)=0$ whenever the coderivation $Q$ is constant or linear, $[s(x),-]$ is a derivation with respect to $\circ$ for all $x\in L$. In particular, 
\[ [s(x),s^\bot(y)\circ s^\bot(z)]=[s(x),s^\bot(y)]\circ s^\bot(z)+(-1)^{|x||y|}s^\bot(y)\circ[s(x),s^\bot(z)]=0\quad\mbox{for all $x\in L$.} \]
%Finally, the coderivation $s^\bot(y)\circ s^\bot(z)$ satisfies 
Since moreover $(s^\bot(y)\circ s^\bot(z))_0(1)=ps^\bot(y)s^\bot(z)(1)=s^\bot(y)_1(z)=(-1)^{|y||z|}z\triangleright y=y\triangleleft z$, by the second part of the proposition $s^\bot(y)\circ s^\bot(z)=s^\bot(y\triangleleft z)$.
\end{proof}
\begin{definition} The \emph{exponential automorphism} $E:SL\to SL$ of the symmetric coalgebra $SL$ (associated to the pre-Lie product $\triangleright$ on $L$) is the composition $E:SL\xrightarrow{\operatorname{sym}}UL\xrightarrow{\eta}SL$, with $\operatorname{sym}$ as in Theorem \ref{PBW} and $\eta$ as in Theorem \ref{preLiePBW}. The inverse $E^{-1}:SL\to SL$ is the \emph{logarithmic automorphism} of $SL$.\end{definition}
The choice of names is due to the fact that, under suitable hypotheses ensuring convergence, for instance, if $\triangleright$ is continuous with respect to a complete filtration on $L$, the induced bijections between group-like elements of $SL$ are the usual (left) pre-Lie exponential and logarithmic maps. Recall \cite{chrono,pLDF} that (under the same suitable hypotheses as before) the (left) \emph{pre-Lie exponential map} $e_{\triangleright}^--1:L^0\to L^0$ is defined by \[ e_\triangleright^--1:L^0\to L^0:a\to e_\triangleright^a-1:=a + \frac{1}{2}a\triangleright a +\frac{1}{6}a\triangleright(a\triangleright a)+\cdots+\frac{1}{n!}a\triangleright(\cdots\triangleright(a\triangleright a)\cdots)+\cdots.\]
This is a bijection with inverse the (left) \emph{pre-Lie logarithm} $\log_\triangleright(-+1):L^0\to L^0$. The latter may be computed recursively via the identity $x:=\log_\triangleright(a+1)=\sum_{n\geq0}\frac{B_n}{n!}(\nabla_{x})^n(a)$, cf. \cite{chrono}, where $B_0=1$, $B_1=-\frac{1}{2}$, $B_2=\frac{1}{6}$, $B_3=0$, $B_4=-\frac{1}{30}$, ..., are the Bernoulli numbers: the first few terms are
\[\log_{\triangleright}(a+1)=a-\frac{1}{2}a\triangleright a+ \frac{1}{4}(a\triangleright a)\triangleright a+ \frac{1}{12}a\triangleright (a\triangleright a)+\cdots.\]  
\begin{remark}\label{rem:fake1} In the previous identities we are not using $1$ in the left hand side to denote a particular element of $L$, just as a formal symbol to remind us that in the case of a unitary associative algebra, seen as a left pre-Lie algebra, we recover the power series expansions of the usual exponential and logarithmic functions $e^a-1$, $\log(a+1)$.
\end{remark}
Given a degree zero element $a\in L^0$, we denote by $a_\odot^n=a\odot\cdots\odot a$ and $a_\ast^n=a\ast\cdots\ast a$ the $n$-th power of $a$ in $SL$ and $UL$ respectively, and by $e_\odot^a=\sum_{n\geq 0}\frac{1}{n!}a_\odot^n$, $e_\ast^a=\sum_{n\geq 0}\frac{1}{n!}a_\ast^n$ the corresponding group-like elements of $SL$ and $UL$ (to be more precise, the latters should be considered as elements in the completed coalgebras $\widehat{SL}$ and $\widehat{UL}$, cf. \cite{RHT}, but with a little abuse we shall overlook this point). We have $\operatorname{sym}(e_\odot^a)=e_\ast^a$ and $\eta(e_\ast^a)=E(e_\odot^a)=e_\odot^{(e_\triangleright^a-1)}$. The first identitity is clear, to prove the second it suffices to show that $p\eta(e_\ast^a)=e_\triangleright^a-1$, where $p:SL\to L$ is the natural projection: we have $p\eta(1)=0$ and an easy induction shows $p\eta(a_\ast^n)=ps(a)^n(1)=\nabla_a(ps(a)^{n-1}(1))=\nabla_a^{n-1}(a)$ for all $n\geq1$, from which the claim follows. 
\begin{remark} By a standard argument (cf., for instance, \cite{Bering}), we see in particular that $E,E^{-1}$ are induced by the pre-Lie exponential and logarithm via the usual polarization trick. More precisely, $E$ is given in Taylor coefficients $e_n:L^{\odot n}\to L$, $n\geq1$, by
\begin{equation}\label{E} e_1=\operatorname{id}_L,\qquad e_n(x_1\odot\cdots\odot x_n)=\frac{1}{n!}\sum_{\sigma\in S_n}\varepsilon(\sigma)x_{\sigma(1)}\triangleright(\cdots\triangleright(x_{\sigma(n-1)}\triangleright x_{\sigma(n)})\cdots),  \end{equation}
while the Taylor coefficients of $E^{-1}$ are determined recursively by $e_1^{-1}=\id_L$, and for $n\geq2$
\begin{multline*}e^{-1}_n(x_1\odot\cdots\odot x_n)=\sum\limits_{k=1}^{n-1}\frac{B_k}{k!}\sum_{i_1+\cdots+i_k=n-1}\sum_{\sigma\in S(i_1,\ldots,i_k,1)}\varepsilon(\sigma)\\e^{-1}_{i_1}(x_{\sigma(1)}\odot\cdots\odot x_{\sigma(i_1)})\triangleright\left(\cdots\triangleright\left(e^{-1}_{i_k}(x_{\sigma(n-i_k)}\odot\cdots\odot x_{\sigma(n-1)})\triangleright x_{\sigma(n)}\right)\cdots\right).\end{multline*}\end{remark}

Following \cite{pLDF}, we shall denote by $(-+1)\circledcirc -:L^0\times L\to L: (a,x)\to (a+1)\circledcirc x \,:=\, ps^\bot(x)(e_\odot^a)=\sum_{n\geq0}\frac{1}{n!}\{ a,\ldots,a;x\}$ and call it the \emph{circle product} on $L$ (as in Remark \ref{rem:fake1}, we treat $1$ as a formal symbol), moreover, we shall denote by $-\bullet-:L^0\times L^0\to L^0$ the Baker-Campbell-Hausdorff product on the Lie algebra $(L^0,[-,-])$. The previous setup implies rather naturally the following computation, from \cite{pLDF}, of the formal group law on $L^0$ associated to the left pre-Lie product on $L$, that is, the transfer of $\bullet$ via the pre-Lie exponential and logarithm.
\begin{theorem}\label{th:formalgrouplaw} For all $a,b\in L^0$ (under suitable hypotheses ensuring convergence) we have
\[e_\triangleright^{\log_\triangleright(a+1)\bullet\log_\triangleright(b+1)}-1\,\,=\,\, a + (a+1)\circledcirc b\,\,=\,\, a + e^{\nabla_{\log_\triangleright(a+1)}}(b).\]
\end{theorem}   
The first identity is \cite[Sec. 4, Th. 2]{pLDF} while the second is \cite[Sec. 4, Prop. 3]{pLDF}, the identity between the left and the right hand side was proved in \cite{chrono}.
\begin{proof} For simplicity, we put $x:=\log_\triangleright(a+1)$ and $y:=\log_\triangleright(b+1)$. We have $e_\triangleright^{x\bullet y}-1=p\eta\left(e_\ast^{x\bullet y}\right)=p\eta(e_\ast^{x}\ast e_\ast^{y})$: since $e_\ast^y=\eta^{-1}(e_\odot^b)$, we see that $e_\triangleright^{x\bullet y}-1=p\eta (e^{x\ast-})\eta^{-1}(e_\odot^b)=pe^{s(x)}(e_\odot^b)$, where the exponentials $e^{x\ast-}, e^{s(x)}$ are taken in $\End(UL)$ and $\End(SL)$ respectively. We have $p(e_\odot^b)=b$, $ps(x)(e_\odot^b)=p\sigma_x(e_\odot^b)+p\nabla_x(e_\odot^b)=x+x\triangleright b$ and by induction \[ps(x)^n(e_\odot^b)=p(\sigma_x+\nabla_x)s(x)^{n-1}(e_\odot^b)=\nabla_x(ps(x)^{n-1}(e_\odot^b))=\nabla_x^{n-1}(x+x\triangleright b)\] for all $n\geq2$, thus $pe^{s(x)}(e_\odot^b)=(e_\triangleright^x-1)+e^{\nabla_x}(b)=a+e^{\nabla_{\log_\triangleright(a+1)}}(b)$. 

Reasoning as before, we see moreover $e_\triangleright^{x\bullet y}-1=p\eta (e^{-\ast y})\eta^{-1}(e_\odot^a)=pe^{s^\bot(y)}(e_\odot^a)$. In the right pre-Lie algebra $(L,\triangleleft)$, we put $y_\triangleleft^n=(\cdots(y\triangleleft y)\triangleleft \cdots)\triangleleft y =y_\triangleright^n$, similarly, given a coderivation $Q\in\Coder(SL)$, we put $Q_\circ^n=(\cdots(Q\circ Q)\circ\cdots)\circ Q$. For all $n\geq1$ we have $pQ_\circ^n=pQ^n:SL\to L$: for $n=1$ this is trivial and for $n=2$ it is true by definition of the Nijenhuis-Richardson product, in general, by induction, $pQ_{\circ}^n=p(Q_{n-1}^\circ\circ Q)=pQ_{\circ}^{n-1}Q=pQ^n$. Finally, by Proposition \ref{prop:symmbraces} 
\begin{multline*} pe^{s^\bot(y)}(e_\odot^a)= a +\sum_{n\geq1}\frac{1}{n!}ps^\bot(y)^n(e_\odot^a)=a +\sum_{n\geq1}\frac{1}{n!}ps^\bot(y)_\circ^n(e_\odot^a)=a+\sum_{n\geq1}\frac{1}{n!}ps^\bot(y_\triangleleft^n)(e_\odot^a)=\\=a+\sum_{n\geq1}\frac{1}{n!}ps^\bot(y_\triangleright^n)(e_\odot^a)=a+ps^\bot(e_\triangleright^y-1)(e_\odot^a)=a+ps^\bot(b)(e_\odot^a)=a+(a+1)\circledcirc(b).\end{multline*}
\end{proof}

\begin{definition}\label{def:kaprbrack} We denote by $\sK:\Coder(SL)\to\Coder(SL):Q\to EQE^{-1}$ the twisting by the exponential automorphism of $SL$: this is an automorphism of the graded Lie algebra $\Coder(SL)$, with inverse $\sK^{-1}:\Coder(SL)\to\Coder(SL):Q\to E^{-1}QE$. Given an endomorphism $f:L\to L$, regarded as a linear coderivation on $SL$, we shall call the Taylor coefficients $\sK(f)_n:L^{\odot n}\to L$, $\sK^{-1}(f)_n:L^{\odot n}\to L$ respectively the \emph{Kapranov brackets} and the \emph{Koszul brackets} on the pre-Lie algebra $L$ associated to $f$.\end{definition}
\begin{remark}\label{rem:koszul} Given a graded commutative algebra $(A,\cdot)$, regarded as a left pre-Lie algebra, and $f:A\rh A$, the $\sK^{-1}(f)_n$ are the usual Koszul brackets on $A$ associated to $f$ \cite{koszul}: this follows directly by results of Markl \cite{Markl1,Markl2}. If we drop graded commutativity of $(A,\cdot)$ but maintain associativity, we recover the non-commutative Koszul brackets considered in \cite{Bering,HDB,Manettikoszul} (actually, as in \cite{koszul} the first two references deal with unitary algebras, and consider slightly different brackets twisted by the unit $1_A\in A$).
\end{remark}
%\begin{remark} The inverse $E^{-1}:SL\rh SL$ is explicitly given in Taylor coefficients $e^{-1}_n:L^{\odot n}\rh L$, $n\geq1$, by $e^{-1}_1=\operatorname{id}_L$ and then for $n\geq2$ by the recursion where the $B_k$ are the Bernoulli numbers. For instance $e^{-1}_2(x_1\odot x_2)=\sum_{\sigma\in S_2}-\frac{\varepsilon(\sigma)}{2}x_{\sigma(1)}\triangleright x_{\sigma(2)}$ and\[e^{-1}_3(x_1\odot x_2\odot x_3)=\sum_{\sigma\in S_3}\varepsilon(\sigma)\left(\frac{1}{4}(x_{\sigma(1)}\triangleright x_{\sigma(2)})\triangleright x_{\sigma(3)}+\frac{1}{12}x_{\sigma(1)}\triangleright(x_{\sigma(2)}\triangleright x_{\sigma(3)})\right).\]Notice how the induced $L^0\rh L^0:x\rh\sum_{n\geq1}\frac{1}{n!}e_n(x^{\odot n})$ and $L^0\rh L^0:x\rh\sum_{n\geq1}\frac{1}{n!}e^{-1}_n(x^{\odot n})$ are respectively the pre-Lie exponential and the pre-Lie logarithm of the left pre-Lie algebra $L^0$: we refer to \cite{esi}, Section 1.1, in particular formulas (8)-(10). It follows from this and a standard polarization argument that the above defined $E^{-1}$ is in fact the inverse to $E$ .  \end{remark}
\begin{proposition}\label{prop:main} Given a graded left pre-Lie algebra $L$ and a derivation $d\in\Der(L,[\cdot,\cdot])$ of the associated graded Lie algebra, the Kapranov brackets $\sK(d)_n:L^{\odot n}\rh L$ are determined by the recursion
\begin{equation}\label{recursivedefinition} \left\{\begin{array}{l} \sK(d)_0 = 0,\,\,\,\sK(d)_1 = d, \\ \sK(d)_2(x\odot y) = \nabla_{dx}(y) - [d,\nabla_x](y), \\ \sK(d)_{n+1}(x\odot y_1\odot\cdots\odot y_n) =  - [\sK(d)_n,\nabla_x](y_1\odot\cdots\odot y_n)\qquad\mbox{for $n\geq2$,}    \end{array}\right. \end{equation}
where the bracket in the right hand side is the Nijenhuis-Richardson bracket.
\end{proposition}
\begin{proof} $EdE^{-1}(1)=Ed(1)=0$, thus $\sK(d)_0=0$. We can write the above recursion in the more compact form
\begin{equation}\label{compactrecursion} [\sK(d),s(x)]=[\sK(d),\sigma_x+\nabla_x] =\sigma_{dx}+\nabla_{dx}=s(dx)\,\,\,\,\,\,\,\,\,\,\forall x\in L.\end{equation}
In fact, taking the induced identity between the $n$-th Taylor coefficients in \eqref{compactrecursion}, $n\ge0$, we recover the recursive definition of $\sK(d)_{n+1}$ in \eqref{recursivedefinition}, cf. Remark \ref{brwithcostcoder}. We have $\sK(d)=\eta(\operatorname{sym} d\operatorname{sym}^{-1})\eta^{-1}$. Since $d$ is a Lie algebra derivation it induces a biderivation of the bialgebra $UL$, which we denote by $d_{UL}$, and it is straightforward to check $\operatorname{sym} d\operatorname{sym}^{-1}=d_{UL}$. Finally,  
\[ [\sK(d),s(x)]=[\eta(d_{UL})\eta^{-1},\eta(x\ast-)\eta^{-1}]=\eta[d_{UL},x\ast-]\eta^{-1}=\eta(dx\ast-)\eta^{-1}=s(dx), \]
proving \eqref{compactrecursion} and therefore the proposition.
\end{proof}
\begin{remark}\label{rem:old} It would be not a priori obvious (if not for the proposition) that the brackets $\sK(d)_n$ defined by \eqref{recursivedefinition} are graded symmetric: in fact, this follows from the hypothesis $d\in\Der(L,[-,-])$. For instance 
\[ \nabla_{dx}(y)-[d,\nabla_x](y) = dx\triangleright y + (-1)^{|x||d|}x\triangleright dy - d(x\triangleright y), \]
in other words, $\nabla_{dx}-[d,\nabla_x]$ (thus $\sK(d)_2$, if $d\in\Der(L,[-,-])$) measures how far is $d$ from satisfying the Leibniz rule with respect to the pre-Lie product $\triangleright$: clearly, this is graded symmetric in $x$ and $y$ if and only if $d$ is a derivation of the associated Lie bracket. 
%We suppose inductively to have shown graded symmetry of $\Phi(d)_i$ for all $2\leq i\leq n$, then $\Phi(d)_{n+1}$ is clearly graded symmetric in the last $n$ arguments, so we have to show graded symmetry in the first two, that is to say, we have to show that \begin{equation}\label{symmetryofbrackets} [[\Phi(d)_n,\nabla_x],\sigma_y]-(-1)^{|x||y|}[[\Phi(d)_n,\nabla_y],\sigma_x]=0\qquad\forall x,y\in L, n\geq2.\end{equation} We assume for simplicity $n\geq 3$, then we compute \begin{multline*} [[\Phi(d)_n,\nabla_x],\sigma_y]-(-1)^{|x||y|}[[\Phi(d)_n,\nabla_y],\sigma_x]=\,\,\mbox{by Jacobi }\,\,=\\=[\Phi(d)_n,\sigma_{\nabla_x(y)}]+(-1)^{|x||y|}[[\Phi(d)_n,\sigma_y],\nabla_x]-(-1)^{|x||y|}[\Phi(d)_n,\sigma_{\nabla_y(x)}]-[[\Phi(d)_n,\sigma_x],\nabla_y]=\\=\,\,\mbox{by \eqref{recursivedefinition} }\,\,=[\Phi(d)_n,\sigma_{[x,y]}]+[[\Phi(d)_{n-1},\nabla_x],\nabla_y]-(-1)^{|x||y|}[[\Phi(d)_{n-1},\nabla_y],\nabla_x]=\\=\,\,\mbox{by \eqref{recursivedefinition} and Jacobi }\,\,=-[\Phi(d)_{n-1},\nabla_{[x,y]}]+[\Phi(d)_{n-1},[\nabla_{x},\nabla_y]]=0\,\,\,\,\,\mbox{by \eqref{preLiecondition2}}.\end{multline*}The computation in the remaining case $n=2$ is similar, only slightly more involved.
\end{remark}
\begin{remark}\label{rem:kapindofmetr}
Given two left pre-Lie products $\triangleright$ and $\blacktriangleright$ on $L$ with the same associated graded Lie bracket, we denote by $\nabla_-,\nabla'_-:L\rh\End(L)$, $\eta,\eta':UL\rh SL$, $\sK,\sK':\Coder(SL)\to\Coder(SL)$ the respective data defined as before. The automorphism $G:=\eta'\eta^{-1}:SL\rh SL$ of the symmetric coalgebra $SL$ satisfies $G\sK(Q)G^{-1}=\sK'(Q),\,\forall Q\in\Coder(SL)$. We have
\begin{multline*}
G(x\odot y_1\odot\cdots\odot y_n)=\eta'\eta^{-1}(x\odot y_1\odot\cdots\odot y_n)=\\=\eta'\left( x\ast\eta^{-1}(y_1\odot\cdots\odot y_n)-\eta^{-1}(\nabla_{x}(y_1\odot\cdots\odot y_n))\right)=\\=(\sigma_{x}+\nabla'_{x})G(y_1\odot\cdots\odot y_n)-G\nabla_{x}(y_1\odot\cdots\odot y_n).
\end{multline*}
Taking the corestriction on both sides, we see that $G$ is given in Taylor coefficients $g_n:L^{\odot n}\rh L$, $n\geq1$, by $g_1=\operatorname{id}_L$, and for $n+1\geq2$ by the recursion
\begin{multline}\label{cambiometrica}
g_{n+1}(x\odot y_1\odot\cdots\odot y_n)=\\=\nabla'_{x}\left(g_{n}(y_1\odot\cdots\odot y_n)\right)-\sum_{i=1}^{n}(-1)^{\sum_{j=1}^{i-1}|x||y_j|}g_{n}(y_1\odot\cdots\odot\nabla_{x}(y_i)\odot\cdots\odot y_n).
\end{multline}
For instance $g_2(x\odot y)=\nabla'_x(y)-\nabla_x(y)=x\blacktriangleright y-x\triangleright y$: this is graded symmetric since by hypothesis $\blacktriangleright$ and $\triangleright$ have the same associated Lie bracket.
\end{remark}

Notice that the previous proposition implies $[s(x)-\sK(\ad_x),s(y)]=0$ for all $y\in L$: by Proposition \ref{prop:symmbraces} $s(x)-\sK(\ad_x)=s^\bot(x)$ for all $x\in L$, and then $[\sK(d),s^\bot(x)]=s^\bot(dx)$ for all $x\in L$. Putting this facts together we see that, given $d\in\Der^1(L,[-,-])$ such that $[d,d]=2d^2=0$, the correspondence 
\[s-s^{\bot}:(L,d,[-,-])\times(L,d,[-,-])\to(\Coder(SL),[\sK(d),-],[-,-]):(x,y)\to s(x)-s^\bot(y)\]
is a morphism of dg Lie algebras. According to Example \ref{ex:extensions}, this morphism classifies an $L_\infty$ extension of base $L\times L$ and fiber the $L_\infty$ algebra structure on $L[-1]$ induced by the Kaprnaov brackets $\sK(d)_n$: the total space of this $L_\infty$ extension has underlying tangent complex naturally isomorphic to $C^*(\Delta_1;L,d)$, the complex of non-degenerate cochains on the $1$-simplex $\Delta_1$ with coefficients in $(L,d)$, and we may similarly identify the projection over the base with the pull-back $C^*(\Delta_1;L,d)\to C^*(\de\Delta_1;L,d)\cong (L,d)\times (L,d)$ of cochains and the inclusion of the fiber with the inclusion $(L[-1],-d)\cong C^*(\Delta_1,\de\Delta_1;L,d)\to C^*(\Delta_1;L,d)$ of relative cochains, where $\de\Delta_1\subset\Delta_1$ is the boundary. Accordingly, we shall denote by 
\begin{equation}\label{ext1} C^*(\Delta_1,\de\Delta_1;L,d)_{p-Lie}\to C^*(\Delta_1;L,d)_{p-Lie}\to C^*(\de\Delta_1;L,d)_{p-Lie}\end{equation}
the $L_\infty$ extension classified by $s-s^\bot:L\times L\to\Coder(SL)$.
\begin{remark} Given a dg associative algebra $(A,d,\cdot)$, seen as a left pre-Lie algebra, \eqref{ext1} is the extension of dg associative algebras obtained by tensorizing the extension 
	\[C^*(\Delta_1,\de\Delta_1;\K)\to C^*(\Delta_1;\K)\to C^*(\de\Delta_1;\K),\]
with the dg algebra structure given by the usual differential and the cup product, with $(A,d,\cdot)$.
\end{remark}

We want to compare the $L_\infty$ extension \eqref{ext1} with another one, essentially introduced by Fiorenza and Manetti \cite{FMcone}, which we shall denote by 
\begin{equation}\label{ext2} C^*(\Delta_1,\de\Delta_1;L,d)_{Lie}\to C^*(\Delta_1;L,d)_{Lie}\to C^*(\de\Delta_1;L,d)_{Lie}\end{equation}
The base is again $C^*(\de\Delta_1;L,d)_{Lie}:=(L,d,[-,-])\times (L,d,[-,-])$, while this time the fiber is $C^*(\Delta_1,\de\Delta_1;L,d)_{Lie}:=(L[-1],-d)$ regarded as an abelian $L_\infty$ algebra. The underlying tangent complex of $C^\ast(\Delta_1;L,d)_{Lie}$ is again naturally isomorphic to $C^*(\Delta_1;L,d)$. We denote the classifying morphism by $\Phi-\Phi^\bot:L\times L\to\Coder(SL):(x,y)\to \Phi(x)-\Phi^\bot(y)$: it is given in Taylor coefficients by  $\Phi(x)_0(1)=x=\Phi^\bot(x)_0(1)$, and for $n\geq1$
\begin{equation}\label{Phi} \Phi(x)_n(x_1\odot\cdots \odot x_n)=\frac{(-1)^nB_n}{n!}\sum_{\sigma\in S_n}\varepsilon(\sigma)[\cdots[x,x_{\sigma(1)}]\cdots,x_{\sigma(n)}],\end{equation}
\begin{equation*} \Phi^\bot(x)_n(x_1\odot\cdots \odot x_n)=\frac{B_n}{n!}\sum_{\sigma\in S_n}\varepsilon(\sigma)[\cdots[x,x_{\sigma(1)}]\cdots,x_{\sigma(n)}].\end{equation*}
The fact that $\Phi,-\Phi^\bot,\Phi-\Phi^\bot$ are morphisms of dg Lie algebras follows from Theorem \ref{th:extensions} and the results from \cite{FMcone} (cf. also \cite[Sec. 3]{HDB}: with the definitions given there, the $L_\infty$ algebra $C^*(\Delta_1;L,d)_{Lie}$ is the mapping cocylinder of the identity $\id_L:L\to L$). Since $\Phi:L\to\Coder(SL)$ is a morphism of graded Lie algebras and $\Phi(x)_0(1)=x$ for all $x\in L$, Remark \ref{rem:PBWdata} shows that $\varphi:UL\to SL:x_1\ast\cdots\ast x_n\to \Phi(x_1)\cdots\Phi(x_n)(1)$ is an isomorphism of coaugmented coalgebras: we claim that $\varphi=\operatorname{sym}^{-1}$. To prove the claim, by a standard polarization argument it suffices to show that $e^{\Phi(x)}(1)=\varphi(e_\ast^x)=\varphi\operatorname{sym}(e_\odot^x)=e_\odot^x$ for all $x\in L^0$: this follows immediately from $\Phi(x)_n(x_\odot^n)=0$ for all $n\geq1$. In particular, we see that $\Phi(x)=\operatorname{sym}^{-1}(x\ast-)\operatorname{sym}$, thus 
\[ E\Phi(x)E^{-1}=\eta(\operatorname{sym}\Phi(x)\operatorname{sym}^{-1})\eta^{-1}=\eta(x\ast-)\eta^{-1}=s(x) \]
for all $x\in L$. Similarly, $\Phi^\bot(y)=\operatorname{sym}^{-1}(-\ast y)\operatorname{sym}$ and then $E\Phi^\bot(y)E^{-1}=s^\bot(y)$ for all $y\in L$. Finally, we may regard the isomorphism $E:(SL,d)\to(SL,\sK(d))=(SL,EdE^{-1})$ of dg coalgebras as an isomorphism of $L_\infty$ algebras $E:C^{\ast}(\Delta_1,\de\Delta_1;L,d)_{Lie}\to C^\ast(\Delta_1,\de\Delta_1;L,d)_{p-Lie}$. We have just proved $E(\Phi-\Phi^\bot)E^{-1}=s-s^\bot$: according to Lemma \ref{lem:extensions} this implies the following result, where we denote by $\widetilde{E}:C^*(\Delta_1;L,d)_{Lie}\to C^*(\Delta_1;L,d)_{p-Lie}$ the $L_\infty$ isomorphism associated to $E$ as in \eqref{eq:Etilde}.
\begin{theorem}\label{th:deligne} The diagram
\[\xymatrix{ C^*(\Delta_1,\de\Delta_1;L,d)_{Lie}\ar[r]\ar[d]_E & C^*(\Delta_1;L,d)_{Lie}\ar[r]\ar[d]_{\widetilde{E}} & C^*(\de\Delta_1;L,d)_{Lie} \ar@{=}[d] \\ C^*(\Delta_1,\de\Delta_1;L,d)_{p-Lie}\ar[r] & C^*(\Delta_1;L,d)_{p-Lie}\ar[r] & C^*(\de\Delta_1;L,d)_{p-Lie}} \]
is an isomorphism of $L_\infty$ extensions. 
\end{theorem}

The interest for the $L_\infty$ extension \eqref{ext2} lies in the fact that a 1-cochain $_x\xrightarrow{\,a\,}_y\,\in C^1(\Delta_1;L)$, where $x,y\in L^1, a\in L^0$, is Maurer-Cartan in the $L_\infty$ algebra $C^*(\Delta_1;L,d)_{Lie}$ if and only if $x,y$ are Maurer-Cartan elements of $(L,d,[-,-])$ and $a\in L^0$ is a Gauge equivalence between them, $e^a\ast_G y = x$ with the notations of Remark \ref{rem:Gauge}: this follows from the computations in \cite[Sec. 7]{FMcone}. In other words, the Maurer-Cartan elements of $C^*(\Delta_1;L,d)_{Lie}$ are in bijective correspondence with the arrows in the Deligne groupoid of the dg Lie algebra $(L,d,[-,-])$: by the previous theorem these are also in bijective correspondence with the Maurer-Cartan elements in $C^*(\Delta_1;L,d)_{p-Lie}$ via the isomorphism
\[ \MC(\widetilde{E}):\MC(C^*(\Delta_1;L,d)_{Lie})\to\MC(C^\ast(\Delta_1;L,d)_{p-Lie})\,:\,_x\xrightarrow{\,a\,}_y\quad\to\quad\,_x\xrightarrow{e_\triangleright^a-1}_y \]
Putting $b:=e_\triangleright^a-1$, the $1$-cochain $_x\xrightarrow{\,b\,}_y$ is Maurer-Cartan in $C^\ast(\Delta_1;L,d)_{p-Lie}$ if and only if $x,y$ are Maurer-Cartan elements of $(L,d,[-,-])$ and $p(s(x)-s^\bot(y)+\sK(d))(e_\odot^b)=0$. We obtain the following result, which is \cite[Prop. 5]{pLDF}.
\begin{corollary} Given Maurer-Cartan elements $x,y$ of $(L,d,[-,-])$ and $a\in L^0$, then (as usual, under suitable hypotheses ensuring convergence) $e^a\ast_G y= x$ if and only if 
\[\sum_{n\geq1}\frac{1}{n!}\sK(d)_n((e_\triangleright^a-1)_\odot^n)=e_\triangleright^a\circledcirc y - x\triangleright e_\triangleright^a,\]
where we put, cf. Remark \ref{rem:fake1}, $e_\triangleright^a\circledcirc y:=((e_\triangleright^a-1)+1)\circledcirc y$ and $x\triangleright e_\triangleright^a:=x+x\triangleright(e_\triangleright^a-1)$.
\end{corollary}
\begin{remark}\label{rem:loopspace} From the point of view of homotopy theory $C^\ast(\Delta_1,\de\Delta_1;L,d)_{p-Lie}$, or in other words, $L[-1]$ with the $L_\infty$ structure given by the Kapranov brackets $\sK(d)_n$, may be regarded as a model of the the based loop space of $(L,d,[-,-])$ in the homotopy category of $L_\infty$ algebras. Then it is clear that it should be homotopy abelian: this is analog to the well known fact in rational homotopy theory that the minimal model of an $H$-space, in particular, a based loop space, has trivial differential. Of course, by our results we have the explicit isomorphism $E$ between $C^\ast(\Delta_1,\de\Delta_1;L,d)_{p-Lie}$ and the abelian $L_\infty$ algebra $C^\ast(\Delta_1,\de\Delta_1;L,d)_{Lie}$.  
\end{remark}
We close this section by considering an example from algebra.
\begin{example} Let $(V,Q)$ be an $L_\infty[1]$-algebra: we regard $\Coder(\overline{SV})$ as a left pre-Lie algebra via the opposite of the Nijenhuis-Richardson product $Q\triangleright R = (-1)^{|Q||R|+1}R\circ Q$, and we consider the dg Lie algebra structure on the associated Lie algebra $(\Coder(\overline{SV}),[Q,-],[-,-])$ controlling the deformations of the $L_\infty[1]$ algebra $V$. We sketch the computation of the Kapranov brackets $\sK([Q,-])_n$, leaving to the reader to fill up the details and the signs in the formulas. We have $\sK([Q,-])_1=[Q,-]$, while $\sK([Q,-])_2$ measures the failure of $[Q,-]$ to satisfy the Leibniz rule with respect to $\circ$, cf. Remark \ref{rem:old}: as in the proof of Proposition \ref{prop:symmbraces}, the latter is given (up to a sign) by the associator $A(Q,-,-)$. Finally, it is not hard to see inductively, using the direct computation of $A(Q,-,-)$ as base of the induction and the recursion \eqref{recursivedefinition} for the inductive step, that for all $R_1=(r_{1,1},\ldots,r_{1,k},\ldots),\ldots R_n=(r_{n,1},\ldots,r_{n,k},\ldots)\in\Coder(SL)$ the coderivation $\sK([Q,-])_n(R_1\odot\cdots\odot R_n)$ is given in Taylor coefficients $\sK([Q,-])_n(R_1\odot\cdots\odot R_n)_N:L^{\odot N}\to L$ by 
\begin{multline*}\sK([Q,-])_n(R_1\odot\cdots\odot R_n)_N(x_1\odot\cdots\odot x_N)=\\=\sum_{i_1+\cdots+i_n+k=N}\sum_{\sigma\in S(i_1,\ldots,i_n,k)}\pm q_{k+n}(r_{1,i_1}(x_{\sigma(1)}\odot\cdots)\odot\cdots\odot r_{n,i_n}(\cdots\odot x_{\sigma(i_1+\cdots+i_n)})\odot \cdots\odot x_{\sigma(N)}).\end{multline*}
Given a dg Lie algebra $(L,d,[-,-])$, we denote by $Q=(q_1,q_2,0,\ldots,0,\ldots)$ the associated $L_\infty[1]$ algebra structure on $L[1]$: as in Example \ref{eq:quilconstr}, this is $q_1(s^{-1}l)=-s^{-1}dl$, $q_2(s^{-1}l\odot s^{-1}m)=(-1)^{|l|}s^{-1}[l,m]$. By the above formulas $\sK([Q,-])_n=0$ for all $n\geq 3$, and the resulting dg Lie algebra structure on the space (cf. Remark \ref{decalalge})
\[ \CE^*(L,L)\,:=\,\Coder(\overline{S}(L[1]))[-1]\,=\,\prod_{n\geq1}\Hom(L[1]^{\odot n},L[1])[-1]\,=\,\prod_{n\geq1}\Hom(L^{\wedge n}, L)[-n]\]
coincides (perhaps up to signs) with the usual dg Lie algebra structure on the Chevalley-Eilenberg complex of $L$ with coefficients in the adjoint representation: in particular, the latter is homotopy abelian (which is expected, since it is the derived center of the
dg Lie algebra $L$\footnote{Notice that homotopy abelianity is only claimed over the base field $\K$, and not over the Chevalley-Eilenberg algebra $\CE^\ast(L,\K).$}).
\end{example}
\section{Kapranov's brackets in K\"ahler geometry}\label{sec-examples}

Let $X$ be a hermitian manifold, we denote by $\sA_X$ the de Rham algebra of complex valued smooth forms on $X$, and by $\sA(T_X)$ the $\sA_X$-module of smooth forms with coefficients in the tangent bundle $T_X$. We denote by $D=\nabla+\debar:\sA^{\ast,\ast}(T_X)\rh\sA^{\ast+1,\ast}(T_X)\oplus\sA^{\ast,\ast+1}(T_X)$ the Chern connection on~$\sA(T_X)$ (that is, the only connection compatible with both the metric and the complex structure on $T_X$). Finally, we denote by $(z^1,\ldots,z^d)$ a system of local holomorphic coordinates on $X$ and by $(\frac{\de}{\de z^1},\ldots,\frac{\de}{\de z^d})$ the corresponding local frame of $T_X$.

Given $\alpha\in\sA^{p,q}(T_X)$, the contraction operator $\bi_\alpha\in\End^{p-1,q}(\sA(T_X))$ is defined as follows: if in local coordinates $\alpha=\sum_i \alpha^i\ten\frac{\de}{\de z_i}$ and $\beta=\sum_j\beta^j\ten\frac{\de}{\de z_j}$, then $\bi_\alpha(\beta) = \sum_j \left(\sum_i\alpha^i\wedge(\frac{\de}{\de z_i} \contr \beta^j)\right)\ten\frac{\de}{\de z_j}$, where we denote by $\contr$ the contraction of forms with vector fields and by $\wedge$ the exterior product of forms. A straightforward computation shows that \begin{equation*}\label{eq:debarvscontr}[\debar,\bi_\alpha]=\bi_{\debar \alpha}.\end{equation*}

Recall that $\sA(T_X)$ carries a natural structure of (bi)graded Lie algebra $(\sA(T_X),[\cdot,\cdot])$ induced by the bracket of vector fields, cf. for instance \cite{ManettiRendiconti}. Given $\alpha\in\sA^{p,q}(T_X)$ we introduce the operators
\begin{equation*}D_\alpha:=[\bi_\alpha,D]\in\End^{p+q}(\sA(T_X))\quad\mbox{and}\qquad\nabla_\alpha:=[\bi_\alpha,\nabla]\in\End^{p,q}(\sA(T_X)).\end{equation*}
If $\bi_\alpha(\beta)=\bi_\beta(\alpha)=0$, in particular for all $\alpha,\beta\in\sA^{0,\ast}(T_X)$, the usual Cartan identities $[\bi_\alpha,\bi_\beta]=0$ and $[D_\alpha,\bi_\beta]=\bi_{[\alpha,\beta]}$ hold, moreover, since $D_\alpha=[\bi_\alpha,\nabla+\debar]=\nabla_\alpha+(-1)^{|\alpha|}\bi_{\debar \alpha}$,
\begin{equation}\label{eq:Cartan} [\nabla_\alpha,\bi_\beta]=\bi_{[\alpha,\beta]}\qquad\forall\alpha,\beta\in\sA^{0,\ast}(T_X).\end{equation}

We show that the bilinear product 
\[ \triangleright:\sA^{0,\ast}(T_X)\ten\sA^{0,\ast}(T_X)\rh \sA^{0,\ast}(T_X):\alpha \ten\beta\rh\alpha\triangleright \beta :=\nabla_\alpha(\beta)  = D_\alpha(\beta),\]
defines a graded left pre-Lie algebra structure on $\sA^{0,\ast}(T_X)$ precisely when the hermitian metric on $X$ is K\"ahler: in this case, the associated Lie bracket is the usual one.

As well known \cite{Koba}, the curvature $D^2\in\End^2(\sA(T_X))$ is $\sA_X$-linear, thus in local coordinates $D^2(\sum_j\beta^j\otimes\frac{\de}{\de z^j})=\sum_i(\sum_j\beta^j\wedge\Omega^i_j)\otimes\frac{\de}{\de z^i}$, where the 2-forms $\Omega^i_j\in\sA^2_X$ are locally defined by $D^2(\frac{\de}{\de z^j})=\sum_i \Omega^i_j\otimes\frac{\de}{\de z^i}$. For the Chern connection of an hermitian manifold we know moreover that the $\Omega^i_j$ are $(1,1)$-forms \cite{Koba}, thus $D^2=\frac{1}{2}[\nabla+\debar,\nabla+\debar]\in\End^{1,1}(\sA(T_X))$, and looking at the bidegrees
\begin{equation}\label{eq:a} D^2=[\debar,\nabla],\qquad 0=[\nabla,\nabla].\end{equation}
By the Jacobi identity $[\nabla_\alpha,\nabla]=[[\bi_\alpha,\nabla],\nabla]=0$, $\forall\alpha\in\sA(T_X)$, and by the Cartan identity \eqref{eq:Cartan}
\[ [\nabla_\alpha,\nabla_\beta] = [\nabla_\alpha, [\bi_\beta,\nabla] ] = [[\nabla_\alpha,\bi_\beta],\nabla] = [\bi_{[\alpha,\beta]},\nabla] = \nabla_{[\alpha,\beta]},\qquad\forall\alpha,\beta\in\sA^{0,\ast}(T_X). \]
On the right hand side we have the bracket $[\alpha,\beta]$ induced by the one of vector fields, the pre-Lie identity \eqref{preLiecondition2} holds if this coincides with the commutator of $\triangleright$ \[[\alpha,\beta]=\nabla_\alpha(\beta)-(-1)^{|\alpha||\beta|}\nabla_\beta(\alpha)=D_\alpha(\beta)-(-1)^{|\alpha||\beta|}D_\beta(\alpha),
\qquad\forall\alpha,\beta\in\sA^{0,\ast}(T_X).\]
In other words, $\triangleright$ is a left pre-Lie product on $\sA^{0,\ast}(T_X)$ if and only if $D$ is torsion free, but as well known \cite{Koba} this is equivalent to the hermitian metric on $X$ being K\"ahler. We assume in the remainder that $X$ is a K\"ahler manifold.

We notice that $\debar\in\Der(\sA^{0,\ast}(T_X),[-,-])$, in fact $(\sA^{0,\ast}(T_X),\debar,[-,-])$ is the Kodaira-Spencer dg Lie algebra controlling the infinitesimal deformation of the complex structure on $X$ \cite{ManettiRendiconti}. We are in the setup of Proposition \ref{prop:main}, so the brackets $\sK(\debar)_n$, defined as in \eqref{recursivedefinition}, induce an $L_\infty$ algebra structure on $\sA^{0,\ast}(T_X)[-1]$: we denote, as in the previous section, this $L_\infty$ algebra by $C^*(\Delta_1,\de\Delta_1;\sA^{0,\ast}(T_X),\debar)_{p-Lie}$.

Next we recall the construction of the $L_\infty$ algebra structure on $\sA^{0,\ast}(T_X)[-1]$ by Kapranov \cite{Ka}. We can form a bundle of cocommutative coalgebras $ST_X=\oplus_{n\geq0}T_X^{\odot n}$ and a bundle of Lie algebras $\Coder(ST_X)$ over $X$  as in Section \ref{Sec-L_inftyalg}, this time working in the symmetric monoidal category $\mathbf{Bnd}_X$ of holomorphic vector bundles over $X$. $\Coder(ST_X)$ is isomorphic to $\prod_{n\geq0}\Hom(T_X^{\odot n}, T_X)$ as a holomorphic vector bundle, where the symmetric powers and the internal $\Hom(-,-)$ are taken in the category $\mathbf{Bnd}_X$. Looking at the Dolbeault complexes, we have
\[ \sA^{0,\ast}(\Coder(ST_X))\cong\prod_{n\geq0}\sA^{0,\ast}(\Hom(T_X^{\odot n}, T_X)).\]
For all $n\geq0$ it is defined a morphism of dg spaces
\begin{equation}\label{eq:psi} \Psi: \sA^{0,\ast}(\Hom(T_X^{\odot n}, T_X))\rh\Hom(\sA^{0,\ast}(T_X)^{\odot n},\sA^{0,\ast}(T_X)).\end{equation}
For $n=0$ both the left and right hand side become $\sA^{0,\ast}(T_X)$ and $\Psi$ is the identity, for $n\geq1$ it sends $R_n\in\sA^{0,\ast}(\Hom(T_X^{\odot n},T_X))$ to the composition
\[\Psi(R_n):\sA^{0,\ast}(T_X)^{\odot n}\xrightarrow{-\otimes R_n}\sA^{0,\ast}(T_X^{\odot n})\otimes\sA^{0,\ast}(\Hom(T_X^{\odot n}, T_X))\xrightarrow{}\sA^{0,\ast}(T_X)\]
induced by the wedge product of forms and the contraction $T_X^{\odot n}\ten \Hom(T_X^{\odot n}, T_X)\rh T_X$.

\begin{remark} We notice that the brackets $\Psi(R_n)$ are $\sA^{0,\ast}_X$-multilinear in the following graded sense:
\begin{equation}\label{eq:gradedlinear} \Psi(R_n)(\alpha_1\odot\cdots\odot(\omega\wedge\alpha_k)\odot\cdots\odot \alpha_n)=(-1)^{|\omega|(|R_n|+\sum_{j=1}^{k-1}|\alpha_j|)}\omega\wedge\Psi(R_n)(\alpha_1\odot\cdots\odot\alpha_n),\end{equation}
for all $\alpha_1,\ldots\alpha_n\in\sA^{0,\ast}(T_X)$, $\omega\in\sA^{0,\ast}_X$.
\end{remark}

Finally, there is a dg Lie algebra structure on $\sA^{0,\ast}(\Coder(ST_X))$ induced by the bundle of Lie algebras structure on $\Coder(ST_X)$, and it is easy to see that the various $\Psi$ as in \eqref{eq:psi} assemble to a morphism of dg Lie algebras
\[\Psi:\left(\sA^{0,\ast}(\Coder(ST_X)),\debar,[\cdot,\cdot]\right)\rh\left(\Coder(S\sA^{0,\ast}(T_X)),[\debar,\cdot],[\cdot,\cdot]\right),\]
where in the right hand side we regard $\debar$ as a linear coderivation on $S\sA^{0,\ast}(T_X)$.

The hermitian metric and the connections $D,\nabla$ on $T_X$ induce a hermitian metric and connections, which we still denote by $D$ and $\nabla$, on the associated bundles $\Hom(T_X^{\otimes n}, T_X)$, $n\geq0$: these are compatible and \[D=\nabla+\debar\in\End^{1,0}(\sA(\Hom(T_X^{\otimes n}, T_X)))\oplus\End^{0,1}(\sA(\Hom(T_X^{\otimes n}, T_X)))\]
is the Chern connection on the hermitian bundle $\Hom(T_X^{\otimes n}, T_X)$, cf. \cite{Koba}. Following \cite{Ka}, we define a hierarchy of tensors $R_n\in\sA^{0,1}(\Hom(T_X^{\otimes n}, T_X))$, $n\geq2$, starting with the curvature
\[ R_2=\Omega=\sum_{i,j}\Omega^i_j\,dz^j\otimes\frac{\de}{\de z^i}\in\sA^{1,1}(\End(T_X))\cong\sA^{0,1}(\Hom(T_X^{\otimes 2}, T_X))\]
and then for $n+1\geq3$ by the recursion
\begin{equation}\label{kapranovrecursion} \qquad R_{n+1}=\nabla(R_n)\in\sA^{1,1}(\Hom(T_X^{\otimes n}, T_X))\cong\sA^{0,1}(\Hom(T_X^{\otimes n+1}, T_X)).\end{equation}
In \cite[Prop. 2.5.6]{Ka} there is shown that the tensors $R_n$ are totally symmetric in their holomorphic covariant indices $R_n\in\sA^{0,1}(\Hom(T_X^{\odot n}, T_X))$, $\forall n\geq2$. Moreover, by the proof of \cite[Th. 2.6]{Ka} \[R=(0,0,R_2,\ldots,R_n,\ldots)\in\prod_{n\geq0}\sA^{0,1}(\Hom(T_X^{\odot n},T_X))=\sA^{0,1}(\Coder(ST_X))\]
is a Maurer-Cartan element of the dg Lie algebra $\left(\sA^{0,\ast}(\Coder(ST_X)),\debar,[\cdot,\cdot]\right)$, that is,
\begin{equation*}\label{MCeq}\debar R +\frac{1}{2}[R,R]=0. \end{equation*}
Finally, as $\Psi$ is a morphism of dg Lie algebras, this implies that $\debar+\Psi(R)$ is an $L_\infty[1]$ structure on $\sA^{0,\ast}(T_X)$, where again we are regarding $\debar$ as a linear coderivation on $S\sA^{0,\ast}(T_X)$: in fact,
\[ \frac{1}{2}[\debar + \Psi(R),\debar + \Psi(R)]=[\debar,\Psi(R)]+\frac{1}{2}[\Psi(R),\Psi(R)]=\Psi\left(\debar R + \frac{1}{2}[R,R]\right)=0.\]
This is the $L_\infty[1]$ algebra structure on $\sA^{0,\ast}(T_X)$ defined in \cite{Ka}.

\begin{theorem}\label{th:main} The two $L_\infty[1]$ algebra structures $\sK(\debar)$ and $\debar+\Psi(R)$ on the Dolbeault complex $(\sA^{0,\ast}(T_X),\debar)$ are the same, that is, Kapranov's $L_\infty$ algebra coincides with the one we denoted by $C^*(\Delta_1,\de\Delta_1;\sA^{0,\ast}(T_X),\debar)_{p-Lie}$ in the previous section. In particular, there is an $L_\infty$ isomorphsm 
	\[ E:(\sA^{0,\ast}(T_X),\debar)\rh (\sA^{0,\ast}(T_X),\debar+\Psi(R)) \] 
\[e_1=\operatorname{id}_{\sA^{0,\ast}(T_X)}, \qquad e_n(\alpha_1\odot\cdots\odot\alpha_n)=\frac{1}{n!}\sum_{\sigma\in S_n}\varepsilon(\sigma)\nabla_{\alpha_{\sigma(1)}}\cdots\nabla_{\alpha_{\sigma(n-1)}}(\alpha_{\sigma(n)}),\]
where in the left hand side we regard $(\sA^{0,\ast}(T_X),\debar)$ as an abelian $L_\infty[1]$ algebra.  \end{theorem}

\begin{proof} We have to show $\debar+\Psi(R)=\sK(\debar):=E\debar E^{-1}$, where the Taylor coefficients $\sK(\debar)_n$ are defined by the recursion \eqref{recursivedefinition}: for $n=2$
\[ \sK(\debar)_2(\alpha\odot \beta) = \nabla_{\debar \alpha}(\beta)-[\debar,\nabla_\alpha](\beta)= [\bi_{\debar \alpha},\nabla](\beta)-[\debar,[\bi_\alpha,\nabla]](\beta)=\]
\[ = (-1)^{|\alpha|}[\bi_\alpha,[\debar,\nabla]](\beta) = (-1)^{|\alpha|}\bi_\alpha D^2(\beta). \]
In local coordinates, if $\alpha=\sum_i\alpha^i\otimes\frac{\de}{\de z^i}$, $\beta=\sum_j\beta^j\otimes\frac{\de}{\de z^j}$,
\[ \sK(\debar)_2(\alpha\odot \beta)=\sum_k\left(\sum_{i,j}(-1)^{|\alpha|}\alpha^i\wedge\left(\frac{\de}{\de z^i}\contr(\beta^j\wedge\Omega^k_j)\right)\right)\otimes\frac{\de}{\de z^k}=\]
\[ =\sum_k\left(\sum_{i,j}(-1)^{|\alpha|+|\beta|}\alpha^i\wedge\beta^j\wedge\left(\frac{\de}{\de z^i}\contr\Omega^k_j\right)\right)\otimes\frac{\de}{\de z^k}=\Psi(R_2)(\alpha\odot\beta).\]
The thesis follows inductively by comparing the recursions \eqref{recursivedefinition} and \eqref{kapranovrecursion}.

For all $n\geq2$ the bracket $\sK(\debar)_n$ is $\sA^{0,\ast}_X$-multilinear in the sense of \eqref{eq:gradedlinear}: for $n=2$ it follows by $\sK(\debar)_2=\Psi(R_2)$, in general, by graded symmetry, it suffices to show $\sA^{0,\ast}_X$-linearity in the first variable, which follows inductively by the recursive definition and \begin{equation}\label{eq:nablaA_Xlinear}
\nabla_{\omega\wedge\alpha}(\beta)=\omega\wedge\nabla_\alpha(\beta),\qquad\forall\alpha,\beta\in\sA^{0,\ast}(T_X), \omega\in\sA^{0,\ast}_X.\end{equation}
This reduces the proof of $\Psi(R_n)(\alpha_1\odot\cdots\odot\alpha_n)=\sK(\debar)_n(\alpha_1\odot\cdots\odot\alpha_n)$  $\forall\alpha_1,\ldots\alpha_n\in\sA^{0,\ast}(T_X)$, $n\geq3$, to the case $\alpha_k=\frac{\de}{\de z^{i_k}}$, $k=1,\ldots,n$. Finally, 
a direct computation in local coordinates, using \eqref{kapranovrecursion}, shows $\Psi(R_n)\left(\frac{\de}{\de z^{i_1}}\odot\cdots\odot\frac{\de}{\de z^{i_n}}\right)=[\nabla_{\frac{\de}{\de z^{i_1}}},\Psi(R_{n-1})]\left(\frac{\de}{\de z^{i_2}}\odot\cdots\odot\frac{\de}{\de z^{i_n}}\right)$ for all $n\geq3$, hence by induction and \eqref{recursivedefinition}
\begin{multline*}\Psi(R_n)\left(\frac{\de}{\de z^{i_1}}\odot\cdots\odot\frac{\de}{\de z^{i_n}}\right)=[\nabla_{\frac{\de}{\de z^{i_1}}},\Psi(R_{n-1})]\left(\frac{\de}{\de z^{i_2}}\odot\cdots\odot\frac{\de}{\de z^{i_n}}\right)=
 \\ = -[\sK(\debar)_{n-1},\nabla_{\frac{\de}{\de z^{i_1}}}]\left(\frac{\de}{\de z^{i_2}}\odot\cdots\odot\frac{\de}{\de z^{i_n}}\right) = \sK(\debar)_n\left(\frac{\de}{\de z^{i_1}}\odot\cdots\odot\frac{\de}{\de z^{i_n}}\right).\end{multline*}
\end{proof}

%\begin{remark} As we remarked in the introduction,  both $(\sA^{0,\ast}(T_X),\debar)$ and $(\sA^{0,\ast}(T_X),\debar+\Psi(R))$ can be viewed as $L_\infty[1]$ algebras over $\sA_X$ in the sense of Costello \cite{Costello}, Definition 2.1.2. In general the two are only isomorphic as $L_\infty[1]$ algebras over $\C$, as the maps $e_n$ in the claim of the previous theorem are only $\C$-multilinear. It is not hard to see that the vanishing of the Atiyah class of $T_X$ is a necessary condition for the existence of an $\sA_X$-linear $L_\infty[1]$ isomorphism from $(\sA^{0,\ast}(T_X),\debar)$ to $(\sA^{0,\ast}(T_X),\debar+\Psi(R))$ .\end{remark}

%The previous theorem implies that up to a $\C$-multilinear $L_\infty[1]$ isomorphism Kapranov's $L_\infty[1]$ algebra structure on $\sA^{0,\ast}(T_X)$ only depends on the complex structure on $X$, and not on the choice of a K\"ahler metric: in fact we can deduce a stronger result. 

\begin{proposition}\label{prop:kapindofmetr}
Kapranov's $L_\infty[1]$ algebra structure on $\sA^{0,\ast}(T_X)$ is independent on the choice of a K\"ahler metric up to an $\sA^{0,\ast}_X$-multilinear $L_\infty[1]$ isomorphism (defined recursively as in \eqref{cambiometrica}).
\end{proposition}
\begin{proof} Given two K\"ahler metrics on $X$, we denote by $D=\nabla+\debar$ and $D'=\nabla'+\debar$ the respective Chern connections and by $\sK(\debar)$ and $\sK'(\debar)$ the associated Kapranov brackets on $\sA^{0,\ast}(T_X)$. There is an $L_\infty[1]$ isomorphism $G:(\sA^{0,\ast}(T_X),\sK(\debar))\rh(\sA^{0,\ast}(T_X),\sK'(\debar))$ defined recursively as in Remark~\ref{rem:kapindofmetr}. The Taylor coefficients $g_n$ are all $\sA^{0,\ast}_X$-multilinear: by graded symmetry it suffices to check $\sA^{0,\ast}_X$-linearity in the first variable, which follows by induction using \eqref{cambiometrica} and \eqref{eq:nablaA_Xlinear}.\end{proof}


\begin{thebibliography}{99}

\bibitem{chrono} A. A. Agra\v{c}ev, R. V. Gamkrelidze,
\emph{Chronological algebras and nonstationary vector fields}, Problems in geometry, Vol. 11 (Russian), Akad. Nauk SSSR, Vsesoyuz. Inst. Nauchn. i Tekhn. Informatsii, Moscow, 1980, pp. 135-176; available (in english) at 
\texttt{ http://people.sissa.it/~agrachev/agrachev\_files/chrono.pdf}.

\bibitem{BM} R. Bandiera, M. Manetti, \emph{On coisotropic deformations of holomorphic submanifolds}; \texttt{arXiv:1301.6000v2 [math.AG]}.

\bibitem{HDB} R. Bandiera, \emph{Non-abelian higher derived brackets}, Journal of Pure and Applied Algebra \textbf{219} (2015), no. 8, 3292-3313; \texttt{arXiv:1304.4097 [math.QA]}.

\bibitem{tesi} R. Bandiera, \emph{Higher Deligne groupoids, derived brackets and deformation problems in holomorphic Poisson geometry}, PhD Thesis, University of Rome \vr La Sapienza", 2015; available at \texttt{http://www1.mat.uniroma1.it/ricerca/dottorato/TESI/ARCHIVIO/bandieraruggero.pdf}.

\bibitem{Bering} K. Bering, \emph{Non-commutative Batalin-Vilkovisky algebras, homotopy Lie algebras and the Courant bracket}, Comm. Mat. Phys. \textbf{274} (2007), 297-341; \texttt{arXiv:hep-th/0603116}.

\bibitem{bott} Bott, R., \emph{Some aspects of invariant theory in differential geometry}, in \emph{Differential Operators
on Manifolds}, 49-145, CIME, Edizione Cremonese, 1975 (reprinted in his \emph{Collected
Papers}, Vol. 3, 357-453, Birkhauser, Boston, 1996).

\bibitem{liemodels} U. Buijs, Y. Félix, A. Murillo, D. Tanr\'{e}, \emph{Lie models of simplicial sets and representability of the Quillen functor}; \texttt{ arXiv:1508.01442 [math.AT]}.

\bibitem{Caldararu1} D. Calaque, A. C\v{a}ld\v{a}raru, J. Tu, \emph{PBW for an inclusion of Lie algebras}, Journal of Algebra \textbf{378} (2013), 64-79; \texttt{arXiv: 1010.0985v2 [math.QA]}.

\bibitem{Caldararu2} D. Calaque, A. C\v{a}ld\v{a}raru, J. Tu, \emph{On the Lie algebroid of a derived self-intersection}, Advances in Mathematics \textbf{262} (2014), 751-783; \texttt{arXiv:1306.5260 [math.AG]}.

\bibitem{cartan} H. Cartan, S. Eilenberg, \emph{Homological algebra}, Princeton, Princeton University Press, 1956.

\bibitem{ChLiv} F. Chapoton, M. Livernet, \emph{Pre-Lie Algebras and the Rooted Trees Operad}, Int. Math.
Res. Notes (2001), no.8, 395-408;  	\texttt{arXiv:math/0002069 [math.QA]}.


\bibitem{Xuatiyah} Z. Chen, M. Sti\'enon, P. Xu, \emph{From Atiyah classes to homotopy Leibniz algebras};
    \texttt{ arXiv:1204.1075 [math.DG]}.

\bibitem{GetzChen} X. Z. Cheng, E. Getzler, \emph{Transferring homotopy commutative algebric structures}; \texttt{arXiv:math/0610912 [math.AT]}.

\bibitem{ChLaz2}
J.~ Chuang, A.~ Lazarev, \emph{Combinatorics and formal geometry of the master equation},  Letters in Math. Phys. \textbf{103}, no. 1 (2013), 79-112;
\texttt{arXiv:1205.5970}.

\bibitem{Costello} K. Costello, \emph{A geometric construction of the Witten genus, II}; \texttt{ arXiv:1112.0816 [math.QA]}.

\bibitem{Deligneletter}  P. Deligne, letter to J. Millson, \emph{April 24, 1986}; scanned copy available on J. Millson webpage.

\bibitem{pLDF} V. Dotsenko, S. Shadrin, B. Vallette, \emph{Pre-Lie deformation theory}; \texttt{ arXiv: 1502.03280 [math.QA]}. 


\bibitem{dupont}  J. L. Dupont, \emph{Curvature and characteristic classes}, Lecture Notes in Mathematics \textbf{640}, Springer-Verlag (1978).

\bibitem{postlie} K. Ebrahimi-Fard, A. Lundervold, H. Munthe-Kaas, \emph{On the Lie enveloping algebra of a post-Lie algebra},  	Journal of Lie Theory \textbf{25} (2015), no. 4, 1139-1165; \texttt{arXiv:1410.6350 [math.NA]}


\bibitem{FMcone} D. Fiorenza, M. Manetti,
\emph{$L_{\infty}$ structures on mapping cones}, Algebra \& Number Theory \textbf{1} (2007), 301-330;
\texttt{arXiv:0601312 [math.QA]}.


\bibitem{GF} Gelfand I. M., Kazhdan D. A., Fuks D. B., \emph{The actions of infinite dimensional Lie algebras}, Funct. Anal. Appl. \textbf{6} (1972), 9–13 (reprinted in {Collected papers of I. M. Gelfand}, Vol. 3, 349-353, Springer-Verlag, 1989).

\bibitem{Ger} M. Gerstenhaber, \emph{The cohomology structure of an associative ring}, Ann. of Math. \textbf{78} (1963), 267-288.

\bibitem{Getzler04} E.~Getzler, \emph{Lie theory for nilpotent $L_{\infty}$-algebras.}, Ann. of Math. \textbf{170}, no. 1 (2009), 271-301; \texttt{arXiv:math/0404003v4}.


\bibitem{GoMil1} W.M. Goldman, J.J. Millson,
\emph{The deformation theory of
representations of fundamental groups of compact k\"{a}hler manifolds}, Publ. Math. I.H.E.S. \textbf{67} (1988), 43-96.

\bibitem{Hennion} B. Hennion, \emph{Tangent Lie algebra of derived Artin stacks}; \texttt{arXiv:1312.3167 [math.AG]}.

\bibitem{hinichdescent} V.~Hinich, \emph{Descent of Deligne groupoids}, Internat. Math. Res. Notices ,  no. 5 (1997), 223-239; \texttt{arXiv:alg-geom/9606010v3}.

\bibitem{hinichdgC} V. Hinich, \emph{DG coalgebras as formal stack}, J. Pure Appl. Algebra \textbf{162} (2001), 209-250; \texttt{arXiv:math/9812034v1 [math.AG]}.


\bibitem{iacMan2} D. Iacono, M. Manetti, \emph{Semiregularity and obstructions of complete intersections}, Adv. in Math. \textbf{235} (2013), 92-125; \texttt{arXiv:1112.0
425v4 [math.AG]}.

\bibitem{Ka} M. Kapranov, \emph{Rozansky-Witten invariants via Atiyah classes}, Compositio Mathematica \textbf{115} (1999), 71-113; \texttt{arXiv:alg-geom/9704009}.

\bibitem{Koba} S. Kobayashi, \emph{Differential geometry of complex vector bundles}, Princeton University Press (1987).

\bibitem{KoSo} M. Kontsevich, Y. Soibelman, \emph{Deformation theory I}, draft of the book, available at \texttt{www.math.ksu.edu/~soibel}.

\bibitem{koszul} J.-L. Koszul,
\emph{Crochet de Schouten-Nijenhuis et cohomologie}, Ast\'erisque, (Numero Hors Serie) (1985) 257-271.

\bibitem{LaSt} T. Lada, J. Stasheff, \emph{Introduction to sh Lie algebras for physicists}, Int. J. Theor. Phys. \textbf{32} (1992), 1087-1104; \texttt{arXiv:hep-th/9209099}.

\bibitem{Xuexponential} C. Laurent-Gengoux, M. Sti\'enon, P. Xu, \emph{Exponential map and $L_\infty$-algebra associated to a Lie pair}, C. R. Acad. Sci. Paris, Ser. I \textbf{350} (2012), 817-821 ; \texttt{arXiv:1211.3478 [math.QA]}.

\bibitem{Xukapdg} C. Laurent-Gengoux, M. Sti\'enon, P. Xu, \emph{Kapranov dg-manifolds and Poincar\'e-Birkhoff-Witt isomorphisms}; \texttt{ arXiv:1408.2903 [math.DG]}.


\bibitem{Lazarev} A.~ Lazarev, \emph{Models for classifying spaces and derived deformation theory}, Proc. London Math. Soc. \textbf{109} (2014), 40-64.; \texttt{arXiv:1209.3866v3}.

\bibitem{LurieDAG} J. Lurie, \emph{Formal moduli problems}, 2011; available at \texttt{http://www.math.harvard.
edu/~lurie/papers/DAG-X.pdf}.


\bibitem{esi} D. Manchon, \emph{A short survey on pre-Lie algebras}, E. Schrödinger Institut Lectures in Math. Phys., Eur. Math. Soc., A. Carey Ed. (2011); \texttt{http://math.univ-bpclermont.fr/~manchon/biblio/ESI-prelie2009.pdf}.

\bibitem{EDF} M. Manetti, \emph{Extended deformation functors}, Int. Math. Res. Not. \textbf{14} (2002), 719-756; {\texttt{arXiv:math.AG/9910071}}.

\bibitem{ManettiRendiconti} M. Manetti,
\emph{Lectures on deformations of complex manifolds}, Rend. Mat.
Appl. \textbf{24}, no. 7 (2004), 1-183;
{\texttt{arXiv:math.AG/0507286}}.

\bibitem{Manettiformality} M. Manetti, \emph{On some formality criteria for DG-Lie algebras}, Journal of Algebra \textbf{438} (2015),  90-118; \texttt{arXiv:1310.3048v2 [math.QA]}.

\bibitem{Manettikoszul} M. Manetti, G. Ricciardi, \emph{Universal formulas for higher antibrackets}; \texttt{ arXiv:1509.09032  [math.QA]}.

\bibitem{Markl1} M. Markl, \emph{On the origin of higher braces and higher-order derivations}; \texttt{arXiv:1309.7744 [math.KT]}.

\bibitem{Markl2} M.  Markl, \emph{Higher  braces  via  formal  (non)commutative  geometry}; \texttt{arXiv:1411.6964 [math.AT]}.

\bibitem{methazambon} R. Mehta, M. Zambon, \emph{$L_\infty$-algebra actions}, Diff. Geom. and its applications \textbf{30} (2012), 576-587; \texttt{arXiv:1202.2607v2 [math.DG]}.

%\bibitem{Mnev} P. Mn\"ev, \emph{Notes on simplicial BF theory}, Moscow Mathematical Journal \textbf{9} (2009), no. 2, 371-410; \texttt{ arXiv:hep-th/0610326}. 
 
\bibitem{oudom-guin} J.-M. Oudom,  D. Guin, \emph{On the Lie enveloping algebra of a pre-Lie algebra}, Journal of K-theory \textbf{2} (2008), 147-167; \texttt{ 	arXiv:math/0404457 [math.QA]}.

\bibitem{PridUDDT} J. Pridham, \emph{Unifying derived deformation theories}, Adv. in Math. \textbf{224} (2010), no.3 , 772-826; \texttt{arXiv:0705.0344v6}.

\bibitem{RHT} D. Quillen, \emph{Rational homotopy theory}, Ann. of Math. \textbf{90} (1969), 205-295.

\bibitem{Schl} M. Schlessinger, \emph{Functors of Artin rings}, Trans. Amer. Math. Soc. \textbf{130} (1968), 208-222.

\bibitem{SchSta} M. Schlessinger, J. D. Stasheff, \emph{Deformation theory and rational homotopy type}; \texttt{ arXiv:1211.1647}.

\bibitem{sullivan} D. Sullivan, \emph{Infinitesimal computations in topology}, Publications math\'ematiques de l'I.H.\'E.S. \textbf{47} (1977), 269-331.

\bibitem{vin} E. B. Vinberg, \emph{The theory of homogeneous convex cones}, Transl. Moscow Math. Soc. \textbf{12} (1963), 340-403.

\bibitem{weibel} C. A. Weibel, \emph{An introduction to homological algebra}, Cambridge Studies in Advanced Mathematics \textbf{38}, Cambridge University Press, Cambridge, (1994).

\bibitem{Yu} S. Yu, \emph{The Dolbeault dga of the formal neighborhood of the diagonal}, Journal of Noncommutative Geometry \textbf{9} (2015), no. 1, 161-184; \texttt{ arXiv:1211.1567 [math.AG]}.



\end{thebibliography}
\end{document}